\documentclass[12pt]{article}
\usepackage{amsmath, amsfonts, amsthm, amssymb, bbm, hyperref,enumerate,color}
\usepackage[latin1]{inputenc}
\usepackage{latexsym,pdfsync}

\allowdisplaybreaks
\theoremstyle{plain}
\newtheorem{theorem}{Theorem}[section]
\newtheorem{proposition}[theorem]{Proposition}
\newtheorem{corollary}[theorem]{Corollary}
\newtheorem{lemma}[theorem]{Lemma}
\newtheorem*{assumptionstar}{Assumption}
\newtheorem{assumption}[theorem]{Assumption}

\theoremstyle{definition}
\newtheorem{remark}[theorem]{Remark}

\def\Econd#1#2{\E_{#1}\left[#2\right]}
\def\Pcond#1#2{\P_{#1}\left[#2\right]}
\def\Protau#1{\P_{\tau}\left[#1\right]}
\def\Dist#1{d_{\Zc}\left(#1\right)}
\newcommand\refAssOZx{{\bf (P)}}
\newcommand\refAssZx{{\bf (Z)}}
\newcommand\refAssOZxi{{\bf (P-{\emph i})}}
\newcommand\refAssOZxii{{\bf (P-{\emph {ii}})}}
\newcommand\refassbornetaux{{\bf (L)}}
\newcommand\refassbornetauxprime{{\bf (L')}}
\newcommand\refassbornetauxsecond{{\bf (L'')}}
\newcommand \lp{L^{(p)}}
\newcommand \lpm{L^{(p-1)}}

\def \be{\begin{eqnarray}}
\def \ee{\end{eqnarray}}
\def \b*{\begin{eqnarray*}}
\def \e*{\end{eqnarray*}}
\def \E{\mathbb{E}}
\def \F{\mathbb{F}}
\def \G{\mathbb{G}}
\def \Q{\mathbb{Q}}
\def \R{\mathbb{R}}
\def \P{\mathbb{P}}
\def \1{{\bf 1}}
\def \ep{\hbox{ }\hfill$\Box$}
\def\Ac{{\cal A}}
\def\Bc{{\cal B}}
\def\Ec{{\cal E}}
\def\Fc{{\cal F}}
\def\Gc{{\cal G}}
\def\Oc{{\cal O}}
\def\Tc{{\cal T}}
\def\Zc{{\cal Z}}
\def\eps{\varepsilon}
\def\vs#1{\vspace{#1mm}}
\def\reff#1{{\rm(\ref{#1})}}
\def\x{\times}
\def\pourtout{\mbox{ for all } }
\def\Esp#1{\mathbb{E}\left[#1\right]}
\def\Pro#1{\mathbb{P}\left[{#1}\right]}

\title{First time to exit of a continuous It\^{o} process:  general moment estimates and  
       ${\rm \bf L}_1$-convergence rate for  discrete time approximations}
\author{Bruno Bouchard\thanks{CEREMADE, Universit\'e Paris Dauphine and CREST-ENSAE, place du Mar\'echal de Lattre de Tassigny,
75775 Paris Cedex 16, France. Email: {\tt bouchard@ceremade. dauphine.fr}}~, 
       Stefan Geiss\thanks{Department of Mathematics, University of Innsbruck,
            A-6020 Innsbruck, Technikerstra\ss e 13/7, Austria. Email: {\tt stefan.geiss@uibk.ac.at} }  ~and~  
       Emmanuel Gobet\thanks{CMAP, Ecole Polytechnique and CNRS, Route de Saclay, 91128 Palaiseau cedex, France. Email: {\tt emmanuel.gobet@polytechnique.edu}}}

\begin{document}
\maketitle

\begin{abstract}     
We establish general moment estimates for the discrete and continuous exit times of a general It\^{o} process in terms of  the distance to the boundary. These estimates serve as  intermediate steps to obtain strong convergence results for   the approximation of a continuous exit time  by a discrete counterpart, computed on a grid. In particular, we prove that the discrete exit time of the Euler scheme of a diffusion  converges in the ${\rm \bf L}_1$ norm with an order $1/2$ with respect to  the mesh size. This rate is optimal.
\end{abstract}

\vspace{10mm}

\noindent{\bf Key words:} Exit time; Strong approximation; Euler scheme.
\vspace{5mm}

\noindent{\bf AMS 2000 subject classifications:}  65Cxx; 60F15.


\section{Introduction}

This paper is motivated by the study of the strong convergence rate of the discrete time approximation of the first exit time $\theta$ of a process $Z$ from a non-empty open subset
$\Oc$. 

The interest for numerical discretization of diffusion processes  dates back to the sixties, see \cite{maru:55,mull:56} and \cite{kloe:plat:95} for general references.  Different approaches can be used to approximate the first exit time of a diffusion. We briefly recall them for  the sake of completeness and to make clear the contribution of this paper.
\vs2

{\rm a.} By the very nature of the problem, space discretization schemes naturally appear. The first version is  based on the Walk On Sphere (WOS) schemes introduced in \cite{mull:56}. In  the Brownian motion case one simulates its  position by the first hitting time of a ball contained in the domain and centered at the starting point: the position  is uniformly distributed on the sphere and thus straightforward to sample. The  sampled point is then used as a new starting point. One repeats the above procedure  until  one gets close enough to the boundary of $\Oc$. For a time-homogeneous diffusion process $X$ the scheme is modified using small balls and an Euler-Maruyama approximation. In \cite{mils:96a,mils:98} strong error estimates on the \emph{exit position} $X_{\theta}$ are proved, assuming in particular that the domain $\Oc$ is convex and that the diffusion coefficient satisfies a uniform ellipticity condition. These results do not include an approximation
of the exit time $\theta$.
Weak approximation results -- i.e. for $\E[\varphi(X_{\theta})]$ with $\varphi$ continuous and bounded -- are established in \cite{mils:97}.
\vs1

{\rm b.} For polygonal domains moving from spheres to spheres may not be suitable because of the corners. One has to replace balls by parallelepipeds (tensor products of intervals). Exit times from parallelepipeds are easy to sample. Faure \cite{faur:92a} was probably the first one who developed these ideas. In \cite{mils:tret:99} these ideas are further analyzed for diffusion processes with time-dependency by exploiting small parallelepipeds. Strong error estimates of the exit position and the exit time are established: the  order of convergence of the exit time approximation is $1-\varepsilon$ with respect to~the space step (for any $0<\varepsilon <1$), i.e. equivalently $\frac 1 2 -\varepsilon$ (for any $0<\varepsilon <1/2$) with respect to~the time step, see \cite[Theorem 8.2]{mils:tret:99}. Here again, convexity of $\Oc$ and strong ellipticity were assumed. Related 
{ simulations}
are discussed in \cite{zein:leja:deac:10}. Extensions to non-small parallelepipeds are investigated  in \cite{deac:leja:06}.\vs1

{\rm c.} To maintain a certain simplicity of the simulation, one can alternatively perform the usual Euler scheme on a grid $\pi$ with deterministic time step $|\pi|$ and stop when it exits $\Oc$. This is a crude approximation,  nevertheless the simplest and quickest to use: this is why it has gained much interest in the applied 
{ probability}
community. It results in an order of weak convergence equal to $\frac{ 1 }{ 2}$ with respect to $|\pi|$, see \cite{gobe:00,gobe:meno:10}. Interestingly,  it is shown in \cite{gobe:meno:07} that this order of weak convergence  remains valid for general It\^{o} processes, far beyond the usual diffusion framework in which one can rely on PDE tools to decompose the error. The strong convergence of the exit time is stated in \cite[Theorem 4.2]{gobe:mair:05} but without speed. Finally, note that different techniques can be used to  speed-up the convergence in the  weak sense: sampling the continuous time exit using diffusion bridge techniques \cite{bald:95,gobe:00,bald:cara:02} (possibly with  local modifications of the boundary \cite{gobe:01d,buch:pete:06} or exponential-time stepping \cite{jans:lyth:05}) or using discrete exit times combined with an inward shifting of the boundary \cite{gobe:meno:10}. To our knowledge, no strong error estimates are available for these schemes.  
\vs3

As a matter of fact, until recently  only little was known about the rate of ${\rm \bf L}_1$ convergence of the discrete exit time of an Euler scheme of a diffusion towards the exit time of the exact diffusion, although  there are important fields where the ${\rm \bf L}_1$ criterion is the only relevant one.   As examples let us mention the approximation of backward stochastic differential equations considered in a domain \cite{bouc:meno:09} and the multi level Monte Carlo methods
 \cite{heinrich2001multilevel,gile:08}.
In \cite[Theorem 3.1]{bouc:meno:09} the authors prove that the convergence rate of the discrete exit time of the Euler scheme is of order $\frac{ 1 }{2 }-\varepsilon$ with respect to~$|\pi|$ (for any $0<\varepsilon<1/2$). Because of the aforementioned  applications the question whether one can take $\varepsilon=0$ in the previous estimate has been raised. Also, their arguments are restricted to  finite time horizons and the question whether they could be extended to an infinite time horizon was open.
\vs2
 
In this paper we answer these questions to the positive: the discrete exit time of an Euler scheme converges at the rate $1/2$ in the ${\rm \bf L}_{1}$ norm, even if the time horizon is unbounded,  see Theorem \ref{thm:summary_for_sde}.  In the same theorem we show that the stopped process converges at the rate $1/4$ in ${\rm \bf L}_{2}$. 
Theorem \ref{thm:summary_for_sde} follows from an abstract  version stated in Theorem \ref{thm: error exit times abstract_new}, which we establish in a non-Markovian setting  in the spirit  of \cite{gobe:meno:07}. As a first step of our analysis we provide general controls on the expected time to exit  in terms of the distance to the boundary, see Theorems \ref{thm: bound on tauxa_new} and \ref{thm: bound on taupi0_new} below. They are established  both for continuous exit times and for discrete exit times, i.e. the latter are restricted to take values on a discrete grid.  Essentially, we only use a mild non-characteristic boundary type condition and a uniform bound on the conditional expected times to exit. The fact that, as opposed to most of the papers quoted above, we analyze situations with unbounded time horizon {in a} ${\rm \bf L}_{\infty}$ sense is delicate because the usual finite-time error estimates, e.g. on Euler schemes, blow up exponentially with respect to the time horizon.

In fact our results allow to address much more general problems than the first exit time 
approximations for Markovian stochastic differential equations.  In terms of applications,
many optimal stopping, impulse control, singular control or optimal monitoring problems 
have solutions given by the hitting times of a  domain $\Oc$ by a state process $Z$, see 
e.g.  \cite{shiryaev2007optimal}, \cite{bensoussan1984impulse}, 
\cite{shrevesonertransaction}, {\cite{NeNi90, Fu11,GL14}}. In practice, the process $Z$ is 
only monitored in discrete time and one needs to know how well these hitting times will 
be approximated by counterparts computed on a finite grid. In terms of modeling, there 
is also an increasing need in non-Markovian or infinite dimensional settings, in which 
there is no clear connection between exit times and PDEs with Dirichlet boundary 
conditions. A typical example is the HJM framework for interest rates, see \cite{HJM92}, but this can more generally refer to path-dependent SDEs, see e.g. \cite{Bu00}, or to stochastic evolution equations on Banach spaces, see e.g. \cite{GyMi05}. 

The variety of possible applications motivates the abstract setting of Section \ref{section:moment} and Section \ref{sec: upper bound abstract setting} in which we provide our general moment and approximation estimates on the first exit time of a process $Z$ from a domain $\Oc$.  This process does not need to be neither Markov, nor finite dimensional, we only impose an It\^{o} dynamic for the distance to the boundary and assume that it satisfies a non-characteristic type boundary condition, see Assumption {\bf (P)}. In this general setting, we prove in particular that 
$$
\Esp{  |\theta -\theta^{ \pi}|}= O(|\pi|^{\frac12})
$$ 
where $\theta$ is the first exit time of $Z$, and $\theta^{\pi}$ is its counterpart computed on a time grid $\pi$, with modulus $|\pi|$, see Theorem \ref{thm: error exit times abstract_new} applied to $Z=X=\bar X$. The result remains true when an extra approximation is made on $Z$ and the corresponding distance process converges in ${\rm \bf L}_{1}$ at a rate $1/2$.
We shall  check our general assumptions in details only for  the application to the first exit time approximation of SDEs, see Section \ref{section:L1}. 

We would like to insist on the fact that, even in the simpler context of a Markovian SDE, the advantage of  the abstract results of Section \ref{section:moment} is that they can be applied simultaneously and without extra effort to the original diffusion process and to its Euler scheme. We are not aware of any specific proof that would simplify and shorten our argumentation 
when using the particular setting of Markovian SDEs. 

\vs2
The paper is organized as follows:
In Section \ref{section:moment} we introduce a general set-up followed by the statement of our quantitative results on the moments of the first time to exit. The proof of the main results, Theorems \ref{thm: bound on tauxa_new} and \ref{thm: bound on taupi0_new}, is split into  several subsections. We first establish general Freidlin type inequalities on moments of exit times, which will be controlled in terms of the probability of sub-harmonic paths in Section \ref{sec: sub harmonic}. Estimates on this probability yield to the proof of Theorem \ref{thm: bound on tauxa_new}, that applies to continuous exit times. A final recursion argument is needed to pass from continuous exit times to discrete exit times, see Section \ref{sec: proof of thm: bound on taupi0}.  The application to the exit time approximation error is discussed in  Section \ref{section:L1}, first in an abstract setting, then for the solution of  a stochastic differential equation whose exit time is estimated by the discrete exit time of its Euler scheme.  
\vs2

Throughout this paper, we let $(\Omega,\Fc,\P)$ be a complete probability space supporting a $d$-dimensional Brownian motion $W$. 
We denote by $\F:=(\Fc_{t})_{t\ge 0}$ the right-continuous completion of the natural filtration induced by $W$.	The symbol $\Tc$ denotes  the set of stopping times that are finite a.s.
We write  $\E_{\tau}$ and $\P_{\tau}$ for the conditional expectation and probability,
respectively, given $\Fc_{\tau}$.  
Inequalities between random variables are usually understood in the a.s.-sense without
mentioning it.
Finally, given a vector $a\in \R^d$ or a matrix $A\in \R^{m\times n}$, the notation $|a|$ and 
$|A|$ stands for the {Euclidean} and the Hilbert-Schmidt norm, respectively.


\section{Moment estimates for continuous and discrete exit times}
\label{section:moment}

The main results of this section are Theorems \ref{thm: bound on tauxa_new} and 
\ref{thm: bound on taupi0_new}. They are the basis to prove Theorem \ref{thm: error exit times abstract_new}, which is the main result of the paper in its abstract form. 
 

\subsection{Assumptions}
\label{subsec: pbm formulation abstract}

Let $(\Zc,d_{\Zc})$ be a metric space equipped with the Borel $\sigma$-algebra generated by 
the open sets.
{  In the following we fix an open set $\Oc$} 
of $\Zc$ with
$$\emptyset \not = \Oc \subsetneq \bar \Oc \subsetneq \Zc,$$
in which $\bar \Oc$ denotes the closure of $\Oc$, 
and let $(Z_{t})_{t\ge 0}$ be a continuous $\F$-adapted $\Zc$-valued process starting in $Z_0\equiv z_0 \in \Oc$.  
\vs2
  
The two main  results of this section concern estimates on the time taken by the process $Z$ to reach the boundary of $\Oc$, 
{  where the corresponding exit time takes values in a set $\pi$ which either coincides with 
 $\R_+$ or equals to a countable subset of $\R_+$, that can be thought to be the discretisation 
 points in time of an approximation scheme.}
Therefore the standing assumption of this section is that either
\begin{enumerate}[(a)]
\item $\pi=\R_{+}$,
\item or 
        $\pi$ consists of a strictly increasing sequence $0=t_0<t_1<t_2 < \cdots$ with
       $\lim_n t_n = \infty$ and $|\pi| = \sup_{n\ge 1} |t_n - t_{n-1}| \le 1$.
\end{enumerate}
In both cases, we set
      \begin{equation}\label{eq:phi:phiplus}
      \phi_{t}  :=    \max\{s\in \pi~:~s\le t\}
      \mbox{ and }
      \phi^{ +}_{t} :=   \min\{s\in \pi~:~s\ge  t\},
      \end{equation}
which are the closest points in $\pi$ to the left and to the right of $t$.

\vs2
Our first assumption concerns the path regularity of the process $Z$.  To simplify the notation, we set 
 \be\label{eq: def gamma wrt dZ}
 \gamma(t,s) := d_{\Zc}(Z_t,Z_s),\;\;t,s\ge 0.
 \ee

\vspace{2mm}
 
\begin{assumptionstar}[\bf Z)\;\;(Regularity of $Z$ along $\pi$]
There {  is} a locally bounded map $\kappa$ $:$ $\R_{+}\x (0,\infty) \mapsto \R_{+}$   such that 
	$$
	\Pcond{\tau}{\sup_{\tau \le t\le \tau+T} \gamma(t, \phi_t \lor \tau)  > \rho}
	\le \kappa(T,\rho)|\pi|
	$$  
	for all $\tau \in \Tc$, $T\ge 0$, and $\rho>0$.    
\end{assumptionstar}	 
\vspace{2mm}

{ Although the condition {\bf (Z)} is - so far - a condition on a single fixed time-net
$\pi$, we require the upper bound in a form of a product $\kappa(T,\rho)|\pi|$. As shown in
Lemma \ref{lem: verif ass Z} below, this is a typical form that is also required in our later
computations.}
Our next set of assumptions concerns the 
{ behaviour of the process $Z$ close to}
the boundary $\partial \Oc$ of $\Oc$.   
\vspace{2mm}

 \begin{assumptionstar}[\bf P)\;\;(Distance process $\delta(Z)$] 
There exist  {$L\geq 1$} and an $L$-Lipschitz function $\delta:\Zc\mapsto \R$ such that $\delta>0$ on $\Oc$, 
$\delta=0$   on {$\partial \Oc$}, and  $\delta<0$ on $\bar \Oc^{c}$. In addition, 
the  process $P:=\delta(Z)$ admits the {It\^o process decomposition} 
\begin{equation}
\label{eq: dyna Px}
P_{t}=P_{0} + \int_{0}^{t} b_{s} ds + \int_{0}^{t} a^{\top}_{s} dW_{s}
\end{equation}
for $t\ge 0$, where
\begin{enumerate}[{\rm (i)}]
\item $(P,b,a)$ is a predictable process with values in $[-L,L]^{d+2}$,
\item there is a fixed $r\in (0, L^{-3}/4)$ and a set $\Omega_r\in \Fc$ of measure one, such
      that $|P_t(\omega)|\vee \gamma(t,\phi_t)(\omega)\le r$ implies that
      $|a_t(\omega)| \ge 1/L$ whenever $\omega\in \Omega_r$ and $t\ge 0$.
\end{enumerate}
\end{assumptionstar}
 
Before {we continue, let us comment on} the latter assumptions.
 \begin{remark}\label{rem: explanation on Z and assumptions} {\rm (a)} The process  $P=\delta(Z)$ measures the algebraic distance of $Z$ to the boundary $\partial \Oc$ in terms of the function $\delta$.
 The existence of a signed distance $\delta$ that is 1-Lipschitz can be checked in various
 settings easily (starting from the usual distance one can check whether for all segments
 $[x,y]=\{ z\in \Zc: d_\Zc(x,z)+d_\Zc(z,y) = d_\Zc(x,y)\}$ with $x\in \Oc$ and $y\in (\bar \Oc)^c$ the intersection $[x,y]\cap \partial\Oc$ is
 non-empty), and it can be   modified outside a suitable neighborhood of $\partial \Oc$ in order to be uniformly bounded.
 \vs2
  
~~{\rm (b)}  The It\^{o} decomposition \reff{eq: dyna Px}
 may implicitly impose additional smoothness assumptions on $\partial \Oc$: for instance, if $Z$ is an $\R^d$-valued It\^o process, then $P$ is  also an It\^o process provided that the domain is $C^2$ with compact boundary, see \cite[Proposition 2.1]{gobe:meno:07}. Hence, the condition {\rm (i)} is not too restrictive.
\vs2
 
~~{\rm (c)}  The coefficients $b$ and $a$ may depend on $\pi$. This will be the case in Section \ref{sec: application to SDE} when our abstract results will be applied to an Euler scheme. 
\vs2
   
~~{\rm (d)}  The condition {\rm (ii)} is a uniform non-characteristic boundary condition.
It ensures that the fluctuation of the paths of $Z$ are not  tangential to the boundary. When  $Z$ solves a SDE with diffusion coefficient $\sigma(\cdot)$, i.e. $a^{\top}_{t}=D\delta(Z_{t}) \sigma(Z_{t})$, see Section \ref{sec: application to SDE}, then the natural non-characteristic boundary condition is 
\be\label{eq: cond NC comment}
|D\delta(z) \sigma(z) |\ge 1/L \;\mbox{ if }\; |\delta(z)|\le r,
\ee 
i.e. $|a_{t}|\ge 1/L$ if $|P_{t}|\le r$.  In the case of an Euler scheme $\bar Z$, see  \reff{eq: def schema euler}, we have $\bar a^{\top}_{t}=D\delta(\bar Z_{t}) \sigma(\bar Z_{\phi_{t}})$ and $\bar P_{t}=\delta(\bar Z_{t})$. The natural condition \reff{eq: cond NC comment}  is no more sufficient to ensure that 
$|\bar a_{t}|\ge 1/L$ if $|\bar P_{t}|\le r$. But, by a continuity argument, it is satisfied if the point $\bar Z_{\phi_{t}}$ at which the diffusion coefficient is evaluated  is not too far from the current position $\bar Z_{t}$, i.e. $\gamma(t,\phi_{t})$ is small as well. See Lemma \ref{lem: verif ass O} below.
\end{remark}

\vs2

Now we can define the main objects of this section: given {$\ell\ge 0$}, $\tau \in \Tc$, and 
an integer $p\ge 1$, we set 
\b*
\theta_{\ell}(\tau)       &:= & \inf\{t\ge \tau~:~{P_{t}}\le \ell\}, \\
\theta^{\pi}_{\ell}(\tau) &:= & \inf\{t\ge \tau~:~ t\in \pi,~P_{t}\le  \ell\}, \\
\Phi^{p}_\ell(\tau)       &:= & \Econd{\tau}{( \theta_{\ell}(\tau)-\tau)^{p} }^{\frac1p}, \\
\Phi^{p,\pi}_\ell(\tau)   &:= & \Econd{\tau}{( \theta^{ \pi }_{\ell}(\tau)-\tau)^{p} }^{\frac1p}. 
\e*

Our aim is to provide pointwise estimates on $\Phi^{1}_{0}(\tau)$ and $\Phi^{1,\pi}_{0}(\tau)$.     
 Our arguments require an additional control on the  first conditional moment of the times to exit. 
  
\begin{assumptionstar}[\bf L)\;\;(Uniform bound on expectations of exit times]  
One has that $\Phi^{1,\pi}_{0}(\tau)\le L$ for all $\tau\in \Tc$.
\end{assumptionstar}

In Assumption {\bf (L)} (similarly in Proposition \ref{proposition:reduction_bmo}
and Lemma \ref{lemma: Ecarre tau_new} below) we keep in mind that 
$\theta_0(\tau) \le \theta_0^\pi(\tau)$. Therefore one has
$\Phi^{1}_{0}(\tau)\le \Phi^{1,\pi}_{0}(\tau)$, so that
$\Phi^{1,\pi}_{0}(\tau)\le L$ automatically implies
$\Phi^{1}_{0}(\tau)\le L$.
It should be emphasized  that Assumption {\bf (L)} concerns the given process $(Z_t)_{t\ge 0}$ 
and distance $\delta$, and therefore the fixed distance process $(P_t)_{t\ge 0}$, and
that the same constant $L\ge 1$ as before is taken for notational simplicity.
We refer to \cite[Chapter III, Lemma 3.1]{frei:85} for sufficient conditions ensuring that the exit times of a stochastic differential equation  have finite moments, that are bounded only in terms of the diameter of the domain, the bounds on the  coefficients of the stochastic differential equation and a partial ellipticity condition. 

In Lemma \ref{lemma: Ecarre tau_new} below, we show that \refassbornetaux\ implies that 
$\theta^\pi_{0}(\tau)-\tau$ has finite exponential moments, uniformly in $\tau \in \Tc$.
We conclude this subsection with some equivalent variants of condition {\bf (L)}. The proof is provided in the Appendix.
\vs3

\begin{proposition}\label{proposition:reduction_bmo}
The condition {\bf (L)} is equivalent to either of the following ones:

\begin{enumerate}
\item [{\bf (L')}]
      There is a $L'\ge 1$ such that, for all $\tau\in \Tc$,
      \[ \Phi^{1,\pi}_{0}(\tau)\le L' \mbox{ a.s.~on } \{P_{\tau} > 0\}. \] 
\item [{\bf (L'')}]
      There exist $c>0$ and $\alpha\in (0,1)$ such that, for all $\tau\in \Tc$,
      \[   \P_\tau[\theta^{ \pi }_{0}(\tau)\geq \tau+c]
         \leq \alpha.
        \] 
\end{enumerate}
\end{proposition}


\subsection{First moment control near the boundary}
   
Now we are in a position to state the main results of this section. We will denote by $\Tc^{\pi}$  
the set of stopping times with values in $\pi$. 
Remember that the following can be applied to situations where $\pi=\R_{+}$, in which case assumption \refAssZx\ is automatically satisfied and the extra term $|\pi|^{\frac12}$ below vanishes. 
 
\begin{theorem}\label{thm: bound on tauxa_new} 
Let the assumptions  \refAssZx , \refAssOZx\ and \refassbornetaux\ be satisfied.
\begin{enumerate}[{\rm (a)}]
\item If $\tau\in \Tc^{\pi}$, then
      \[ \Phi^{1}_{0}(\tau)
         \le c_{\reff{thm: bound on tauxa_new}} \Big [ P_{\tau} + |\pi| \Big ] 
             {  \1_{\{P_\tau\ge 0\}}},\]
      where $c_{\reff{thm: bound on tauxa_new}} = c_{\reff{thm: bound on tauxa_new}}(r,L,d,\kappa)>0$.
\item If   $\tau\in \Tc$, then 
      \[ \Phi^{1}_{0}(\tau)
         \le d_{\reff{thm: bound on tauxa_new}}  \Big [ P_{\tau} + |\pi|^\frac{1}{2} \Big ] 
             {\1_{\{P_\tau\ge 0\}}}
             , \]
      where
      $d_{\reff{thm: bound on tauxa_new}} = d_{\reff{thm: bound on tauxa_new}}(r,L,d,\kappa)>0$.
\end{enumerate}
\end{theorem}
\vs3

The proof of this theorem will be given in 
Section \ref{sec: proof of thm: bound on tauxa} below. Its counterpart for discrete exit 
times corresponds to the following statement when $\pi \ne \R_{+}$, and is proved in Section \ref{sec: proof of thm: bound on taupi0}.
\vs3

\begin{theorem}\label{thm: bound on taupi0_new} 
Let the assumptions  \refAssZx , \refAssOZx\ and \refassbornetaux\ be satisfied.
Then there exists an 
$\eps_{\reff{thm: bound on taupi0_new}} = \eps_{\reff{thm: bound on taupi0_new}}(r,L,d,\kappa)>0$ such that 
if $|\pi|\le \eps_{\reff{thm: bound on taupi0_new}}$ then one has 
{\b*
     \Phi^{1,\pi}_{0}(\tau) 
\le  d_{\reff{thm: bound on taupi0_new} }\left[ |P_{\tau}| +|\pi|^{\frac12}  \right] 
     \quad\mbox{for}\quad \tau \in \Tc,
\e* 
where $d_{\reff{thm: bound on taupi0_new} }=d_{\reff{thm: bound on taupi0_new} }(r,L,d,\kappa)>0$.}
\end{theorem}
\vs2

Theorem \ref{thm: bound on tauxa_new} is similar to \cite[Lemma 4.2]{gobe:meno:07}, in which the time horizon is bounded and  the counterpart of  \refAssOZxii\  does not
require $\gamma(\cdot,\phi)\le r$. Our additional requirement yields to a
weaker assumption and explains the presence of the additional  $|\pi|$-terms in our result. We also refer to   \cite[Chapter III, Section 3.3]{frei:85} who considers a Markovian setting
for a uniformly fast exit of a diffusion from a domain.

Theorem \ref{thm: bound on taupi0_new} is of similar nature but is much more delicate to establish.  An attempt to obtain such a result for the Euler scheme of stochastic differential equations on a finite time horizon  can be found in  \cite{bouc:meno:09} by a combination of their Lemmas 5.1, 5.2 and 5.3. However, they were only able to achieve a bound in $O_{{|\pi|\to 0}}(|\pi|^{\frac12-\eps})$ for all $0<\eps<1/2$. We shall comment on this in Section \ref{section:L1} below. {The absolute values on $P_{\tau}$ account for the case where $Z_\tau$ is outside $\Oc$ and $\tau\notin \Tc^{\pi}$ yielding a positive time to exit.}
\vs2
 
The proofs of the above theorems are divided in several steps and provided in the next subsections (see Sections \ref{sec: proof of thm: bound on tauxa} and \ref{sec: proof of thm: bound on taupi0} for the final arguments).  Both     start with  arguments  inspired by \cite{frei:85} and that were already exploited in \cite{bouc:meno:09}. One important novelty is our set of assumptions where we do not use any Markovian hypothesis and where we only assume that the delay to exit is  uniformly bounded  in expectation with respect to the initial time. 
Furthermore, we also refine many important estimates  of  \cite{bouc:meno:09} and use a  new final recursion argument which is presented in Section \ref{sec: proof of thm: bound on taupi0}. This  recursion is crucial in order to recover the bound $O_{{|\pi|\to 0}}(|\pi|^{\frac12})$,
in contrast to the  bound $O_{{|\pi|\to 0}}(|\pi|^{\frac12 - \eps})$ in  \cite{bouc:meno:09}.  

\begin{remark}
Lemma  \ref{lemma: Ecarre tau_new} below implies the same estimates for 
$(\Phi^{p,\pi}_{0}(\tau))^p$ and $(\Phi^{p,\pi}_{0}(\tau))^p$, $p\geq 2$, as obtained for 
$\Phi^{1,\pi}_{0}(\tau)$ and $\Phi^{1,\pi}_{0}(\tau)$ in Theorems \ref{thm: bound on tauxa_new} 
and \ref{thm: bound on taupi0_new}.
\end{remark} 

 \begin{remark}
Theorems \ref{thm: bound on tauxa_new} 
and \ref{thm: bound on taupi0_new}
extend to the case where $\Oc$ is the intersection of 
countable many $(\Oc_{i})_{i\in I}$ satisfying  the assumptions  \refAssOZx\ and \refassbornetaux\ for some family of processes $(P^{i})_{i\in I}$ with the same $L\ge 1$ and 
$r\in (0,L^{-3}/4)$. Indeed, denote by $\Phi^{1,\pi}_{0i}$ and $\Phi^{1}_{0i}$ the counterparts of $\Phi^{1,\pi}_{0}$ and $\Phi^{1}_{0}$ associated to $\Oc_{i}$, $i\in I$, then we have, a.s., 
$$
\Phi^{1}_{0}(\tau)\le \inf_{i\in I} \Phi^{1}_{0i}(\tau) \; \mbox{ and }\;  \Phi^{1,\pi}_{0}(\tau)\le  \inf_{i\in I} \Phi^{1,\pi}_{0i}(\tau)
$$
whenever $\Oc=\cap_{i\in I} \Oc_{i}$.
\end{remark} 

\begin{remark}
Take $d=1$, $\Oc={(-\infty,1)} \subset \Zc=\R$, $\pi=\R_{+}$,  and let $Z=|W|^{2}+z_{0}$
with $1/2<z_{0}< 1$. As distance function take an appropriate $\delta\in C^\infty(\R)$ with 
$\delta$ constant outside $(0,2)$ and $\delta(z)=1-z$ on $[1/2,3/2]$. Then the 
conditions  \refAssZx , \refAssOZx\ and \refassbornetaux\ are satisfied and 
$\Phi^{1}_{0}(0)= \Esp{\theta_{0}(0)}$ $= $ $\Esp{|W_{\theta_{0}(0)}|^{2}}$ $= $ $1-z_{0}=P_0$, which coincides with the upper-bound of Theorem \ref{thm: bound on tauxa_new} up to a multiplicative constant. 
\end{remark} 


\subsection{Freidlin type inequalities on moments of exit times}
 
We start with a-priori estimates inspired by the proof of the exponential fast exit of Freidlin \cite[Lemma 3.3, Chapter 3]{frei:85}:  a uniform bound on the conditional expected times to exit    implies the existence of uniform conditional exponential moments for these  exit times. We adapt 
Freidlin's arguments to our setting.    
 
\begin{lemma}\label{lemma: Ecarre tau_new} 
Let  assumption \refassbornetaux\ hold, $p\ge 1$ be an integer, 
$\lp:=p! L^p$, and $\tau \in \Tc$. Then  we have
\b*
   (\Phi^{p}_{0}(\tau))^{p} 
\le c_{p,\reff{lemma: Ecarre tau_new}} \Phi^{1}_{0}(\tau) &\mbox{ and } &
   (\Phi^{p,\pi}_{0}(\tau))^{p} 
\le c_{p,\reff{lemma: Ecarre tau_new}} \Phi^{1,\pi}_{0}(\tau)
\e*
with $c_{p,\reff{lemma: Ecarre tau_new}}:=p \lpm$.  Consequently,
\b*
       {(\Phi^{p,\pi}_{0}(\tau))^{p}}
&\leq&  L^{(p)},\\
 \E_\tau\left[e^{c(\theta^\pi_{0}(\tau)-\tau)}\right]
&\leq& (1-cL)^{-1},
\e*   
where $c\in [0,L^{-1})$.
\end{lemma}

\proof 
1. The estimates for $\Phi^{p}_{0}(\tau)$ and $\Phi^{p,\pi}_{0}(\tau)$ are obtained in the same way, 
we only detail the  second one by an induction over $p$. The case $p=1$ is an identity.
Assume that the statement is proven for some $p\ge 1$. Observe that, on 
$\{\theta^\pi_{0}(\tau)>t\geq \tau\}=\{\forall s\in [\tau,t]\cap \pi: {Z_s \in \Oc} 
\}$, we have 
 \begin{multline*}
  \theta^\pi_{0}(\tau)
= \inf\{s\ge \tau: s\in \pi,  {Z_s \notin \Oc} \}
= \inf\{s\ge t\lor \tau: s\in \pi,  {Z_s \notin \Oc} \} \\
= \theta^\pi_{0}(t\lor \tau).
\end{multline*}
Hence, for $A\in \Fc_{\tau}$ we can write
\b*
&   & \frac{\Esp{( \Phi^{p+1,\pi}_{0}(\tau))^{p+1}\1_{A}}}{p+1} \\ 
& = & \int_{0}^{\infty }\E \left[\1_{A}(\theta^\pi_{0}(\tau)-t)^{p}
      \1_{\{\theta^\pi_{0}(\tau)>t\ge \tau\}}\right] dt \\
& = & \int_{0}^{\infty }\E \left[\1_{A}\E_{t\lor \tau}[(\theta^\pi_{0}(t\lor\tau)-t\lor \tau)^{p}]
      \1_{\{\theta^\pi_{0}(\tau)>t\ge \tau\}}\right] dt \\
&\le& p! L^{p-1} \int_{0}^{\infty }
      \E \left[\1_{A}\E_{t\lor \tau}[\theta^\pi_{0}(t\lor\tau)-t\lor \tau]
      \1_{\{\theta^\pi_{0}(\tau)>t\ge \tau\}}\right] dt  \\
&\le& p! L^p \int_{0}^{\infty }
      \E \left[\1_{A}
      \1_{\{\theta^\pi_{0}(\tau)>t\ge \tau\}}\right] dt  \\
&\le& \lp\E\left[\1_{A}\E_\tau \left[ \theta^\pi_{0}(\tau)-\tau\right] \right],
\e*
so that the proof is complete because $A\in \Fc_{\tau}$ was arbitrary.
\vs1

2. The consequently part is now obvious.
\ep
  
 
\subsection{An a-priori control in terms of the probability of strictly {sub-harmonic} paths}\label{sec: sub harmonic}
 
Now we provide a control on  $\Phi^{1}_{0}(\tau) $ in terms of the conditional probability 
of
$$
\Ac^\tau_{0}:=\{2Pb +|a|^2 \ge L^{-2}/2  \mbox{ on } [\tau,\theta_{0}(\tau)]\}^{c}.
$$
Intuitively we can say, the more non-degenerate the process $P_t^2$ from  
$\tau$ to $\theta_0(\tau)$ is, the smaller is the time of exit.

\begin{lemma}\label{lem: first bound proof thm bound tauxa} 
Let  assumptions  \refassbornetaux\ and \refAssOZxi \; be satisfied.  
Then there exists a constant
$ c_{\reff{lem: first bound proof thm bound tauxa}}
 =c_{\reff{lem: first bound proof thm bound tauxa}}({L,d})>0$  such that, for all $\tau \in \Tc$,
\b*
     \Phi^{1}_{0}(\tau) 
\le  c_{\reff{lem: first bound proof thm bound tauxa}} \Pcond{\tau}{\Ac^\tau_{0}}.
\e*
\end{lemma}

\proof 
Let $E:= \{ P_\tau \ge 0\}\in \Fc_\tau$ so that $P_{\theta_{0}(\tau)} = 0$ on $E$
and $\Phi^{1}_{0}(\tau) = 0$ on $E^c$. Moreover, on $E$ we obtain that
\b*
       \theta_{0}(\tau)-\tau
&\le& \1_{(\Ac^\tau_{0})^{c}}  {2L^{ 2}} \int_{\tau}^{\theta_{0}(\tau)} (2P_{s}b_{s} +|a_{s}|^2)ds +( \theta_{0}(\tau)-\tau) \1_{\Ac^\tau_{0}}
\\
&=& \1_{(\Ac^\tau_{0})^{c}} 2L^2(|P_{\theta_{0}(\tau)}|^{2}-|P_{\tau}|^{2}) - \1_{(\Ac^\tau_{0})^{c}} 2L^2 \int_{\tau}^{\theta_{0}(\tau)} 2P_{s}a^{\top}_{s} dW_{s}  \\
&&+ ( \theta_{0}(\tau)-\tau) \1_{ \Ac^\tau_{0} }\\
&\le&    - \1_{(\Ac^\tau_{0})^{c}}  4L^2\int_{\tau}^{\theta_{0}(\tau)}  P_{s}a^{\top}_{s} dW_{s} + ( \theta_{0}(\tau)-\tau) \1_{\Ac^\tau_{0}}.
\e* 
Using the bound on $\Phi^{1}_{0}(\tau)$ from assumption \refassbornetaux\  and the bounds from 
assumption \refAssOZxi, we obtain  
$\E \int_0^\infty \1_{\{\tau<s\le \theta_0(\tau)\}} P_s^2 |a_s|^2 ds <\infty$
and, on $E$,
\b*
      \Econd{\tau}{- \1_{(\Ac^\tau_{0})^{c}}  \int_{\tau}^{\theta_{0}(\tau)}  P_{s} a^{\top}_{s} 
      dW_{s}}
& = & \Econd{\tau}{ \1_{\Ac^\tau_{0}}  \int_{\tau}^{\theta_{0}(\tau)} P_{s} a^{\top}_{s} dW_{s}} \\
&\le& L^2\sqrt d \;\P_{\tau}[\Ac^\tau_{0}]^{\frac12}(\Phi^{1}_{0}(\tau))^{\frac12}.
\e*
On the other hand,   Lemma \ref{lemma: Ecarre tau_new} implies
 $$
\Econd{\tau}{(\theta_{0}(\tau)-\tau) \1_{\Ac^\tau_{0}}}\le \Phi^{2}_{0}(\tau)\P_{\tau}[\Ac^\tau_{0}]^{\frac12}\le \left[c_{{2,}\reff{lemma: Ecarre tau_new}} 
{\Phi^{1}_{0}(\tau)}
\P_{\tau}[{\Ac^\tau_{0}}]\right]^{\frac12}.
$$
Combining the above estimates and using the inequality $ab\leq a^2+\frac 14 b^2$ gives,
on $E$,
\b*
       \Phi^{1}_{0}(\tau) 
&\leq& 4 L^4\sqrt d\;\P_{\tau}[{\Ac^\tau_{0}}]^{\frac12}(\Phi^{1}_{0}(\tau))^{\frac12}+
       \left[c_{{2,}\reff{lemma: Ecarre tau_new}} {\Phi^{1}_{0}(\tau)}\P_{\tau}[{\Ac^\tau_{0}}]
       \right]^{\frac12}\\
&\leq& 16 L^8 d\;\P_{\tau}[{\Ac^\tau_{0}}]+\frac 14 \Phi^{1}_{0}(\tau)+c_{{2,}
       \reff{lemma: Ecarre tau_new}}\P_{\tau}[{\Ac^\tau_{0}}]+\frac 14 \Phi^{1}_{0}(\tau),
\e*
which leads to the required result.
\ep 


\subsection{Proof of Theorem \ref{thm: bound on tauxa_new}}
\label{sec: proof of thm: bound on tauxa}

We start by two lemmas before we turn to the proof of Theorem \ref{thm: bound on tauxa_new}.

\begin{lemma}\label{lemma:A_to_B}
Let $\Psi\in \{ \Phi_0^1,\Phi_0^{1,\pi} \}$ and assume that there is a constant $c>0$ 
such that for all $\tau\in \Tc^\pi$ one has that
$$ \Psi(\tau) \le    c \left[ P_{\tau} +|\pi|^{\frac12}  \right]\1_{\{0\le  P_\tau\le   r\}}
             +  L \1_{\{   r<P_\tau  \}}. $$
Then for all $0<\tilde r<r$ there is a $d_{\reff{lemma:A_to_B}}=d_{\reff{lemma:A_to_B}}(r-\tilde r,L,d,c)>0$ such that for all $\tau\in \Tc$ 
one has that
$$ \Psi(\tau) \le d_{\reff{lemma:A_to_B}}
 \left[ |P_{\tau}| +|\pi|^{\frac12}  \right] \1_{\{| P_\tau|\le \tilde r\}} +  L \1_{\{ \tilde r<|P_\tau|  \}}.
$$
\end{lemma}

\proof
The case $\pi=\R_{+}$ is trivial because $\Psi(\tau) \le L$ so that we can assume that 
$\pi\not = \R^+$. Using
$$ \Psi (\tau) \le \Econd{\tau}{\Psi(\phi_\tau^+) + |\pi|} $$ 
and
$$     \Econd{\tau}{|P_{\phi_\tau^{+}} - P_\tau|}
   \le L[1+\sqrt{d}] |\pi|^\frac{1}{2} =: A |\pi|^\frac{1}{2}, $$
{we can conclude by
\b*
&   & \Psi(\tau) \\
&\le& \Econd{\tau}{\Psi(\phi_\tau^+) +|\pi|} 
      \1_{\{ |P_\tau| \le \tilde r\}} + L \1_{\{ \tilde r < |P_\tau|\}} \\
&\le& \Econd{\tau}{c \left[ P_{\phi_\tau^+} +|\pi|^{\frac12}  \right]
      \1_{\{0\le  P_{\phi_\tau^+}\le   r\}}
      +  L \1_{\{   r<P_{\phi_\tau^+}  \}}  +|\pi| } \1_{\{ |P_\tau| \le \tilde r\}} \\
&   & \hspace*{22em}
      + L \1_{\{ \tilde r < |P_\tau|\}} \\
&\le& \left [ c |P_\tau| + [c(1+A)+1] |\pi|^\frac{1}{2} \right ] \1_{\{ |P_\tau| \le \tilde r\}} \\
&   & \hspace*{3em} 
     + L\Pcond{\tau}{r < P_{\phi_\tau^+}, |P_\tau| \le \tilde r}\1_{\{ |P_\tau| \le \tilde r\}}
         + L     \1_{\{ \tilde r < |P_\tau|\}}  \\
&\le& \left [ c |P_\tau| + [c(1+A)+1] |\pi|^\frac{1}{2} \right ] \1_{\{  |P_\tau| \le \tilde r\}} \\
&   & \hspace*{3em}
       + L \Pcond{\tau}{|P_{\phi_\tau^+}-P_\tau| \ge r-\tilde r}
\1_{\{ |P_\tau| \le \tilde r\}}
         + L     \1_{\{ \tilde r < |P_\tau|\}}  \\
&\le& \left [ c |P_\tau| + \left [c(1+A)+1+ \frac{L A}{r-\tilde r}\right ] |\pi|^\frac{1}{2} \right ]
\1_{\{ |P_\tau| \le \tilde r\}}
         + L     \1_{\{ \tilde r < |P_\tau|\}}.
\e*
}\ep
\vs2

Next we control the quantity $\Pcond{\tau}{\Ac^\tau_{0}}$ to make Lemma 
\ref{lem: first bound proof thm bound tauxa}  applicable:
\vs2

\begin{lemma}\label{lem: bound P Axacomp}  
Assume that \refAssZx\ and \refAssOZx\ hold. Then for all $c>0$ there exists an 
$\eta(c)=\eta(c,r,L,d)>0$ such that
\be\label{eq: upper bound PtauActau}
    \P_{\tau}[{\Ac^\tau_{0}}] 
\le \eta(c) P_{\tau} + c\ \Phi^{1}_{0}(\tau) + \kappa\left (\frac{2}{c},r\right ) |\pi|
\;\;\mbox{a.s. on}\;\; \{ P_\tau \in [0,r] \},
\ee
where $\tau \in \Tc^\pi$ and 
$
\Ac^\tau_{0}:=\{2Pb +|a|^2\ge L^{-2}/2  \mbox{ on } [\tau,\theta_{0}(\tau)]\}^{c}.
$
\end{lemma}

\proof
Let $\tilde\theta_r(\tau):=\inf \{ t\ge \tau : P_t = r \} \in [0,\infty]$.
Assumption  \refAssOZxii\ implies $2Pb+|a|^2 \ge  {L^{-2}}/2$ $\P$-a.s.  
on $\{|P|\vee \gamma(\cdot,\phi)\le r\}$ for $r\le L^{-3}/4$.
It follows from the restriction $\tau\in \Tc^\pi$ that on
$$ E:= \{ P_\tau\in [0,r] \} $$
we have, $\P$-a.s., that
\begin{eqnarray*}
            (\Ac^{\tau}_{0})^{c}
&\supseteq& \left \{\sup_{\tau \le t\le \theta_{0}(\tau) }|P_t|\le r \right \}\cap 
            \left \{\sup_{\tau \le  t\le \theta_{0}(\tau) }\gamma(t,\phi_{t}\lor \tau) \le r \right \} 
\\
&\supseteq& \{\theta_{0}(\tau)\le \tilde  \theta_{r}(\tau)\} \cap  
          \left \{\sup_{\tau \le  t\le \theta_{0}(\tau) }\gamma(t,\phi_{t}\lor \tau) \le r \right \}.
\e*
Setting $\Bc_T:=\{\sup_{\tau \le  t\le \tau+T }\gamma (t,\phi_{t}\lor \tau)\le r \}$
for $T:=2c^{- 1 }$, we continue on $E$ with
\begin{eqnarray}
      \P_{\tau}[\Ac^{\tau}_{0}]
&\le& \P_{\tau}[\Ac^{\tau}_{0},\ \theta_{0}(\tau)\leq \tau + T, \Bc_T]+
      \P_{\tau}[\theta_{0}(\tau)>\tau + T]+\P_{\tau}[ \Bc_T^c] \nonumber \\
&\le& \P_{\tau}[\tilde  \theta_{r}(\tau)<\theta_{0}(\tau)\leq \tau + T, \Bc_T]+
      \P_{\tau}[\theta_{0}(\tau)>\tau + T]+\P_{\tau}[ \Bc_T^c] \nonumber \\
&\le& \P_{\tau}[\tilde  \theta_{r}(\tau)<\theta_{0}(\tau)\leq \tau + T, \Bc_T]+
      \frac{c  }{2 } \ \Phi^{1}_{0}(\tau) + \kappa\left (\frac{2}{c},r\right )|\pi|,
      \label{eq: decompo P Ac}
\end{eqnarray}
where the last inequality follows from Chebyshev's inequality and assumption  \refAssZx . 
To treat the first term in \reff{eq: decompo P Ac} we set, {for $T\ge 0$},
$$
\theta_{0,r}^T:=\theta_{0}(\tau) \wedge \tilde \theta_{r}(\tau)\wedge (\tau+T).
$$
In view of assumption \refAssOZx\ we  can define   $\Q\sim \P$ by the density 
	\b*
	\frac{d\Q}{d\P}=H:=\Ec\left( 
	- \int_{\tau}^{.}  \lambda^{\top}_{s} 
	 dW_s
	\right)_{\theta_{0,r}^T} \;,
	\e*
where 
$$
\lambda:=a\ [|a|^{-2}\land L^2]\ b\ \1_{[\tau,\theta_{0,r}^T]}
\;\;\; \mbox{ so that } \;\;\; |\lambda|\leq L^4 \sqrt d=:\lambda_\infty,
$$
and deduce from Girsanov's Theorem (cf. \cite[p.163]{bicht:10}) that 
	$$
	W^{\Q}:=W+ 
	\1_{[\tau,\infty)}\int_{\tau}^{\theta_{0,r}^T\wedge \cdot}   \lambda_{s}ds
	$$ 
is a   Brownian motion associated to $\Q$. For any given $\ell>1$ we obtain
\begin{eqnarray}
&    & \P_{\tau}[\tilde  \theta_{r}(\tau)<\theta_{0}(\tau)\leq \tau + T, \Bc_T]  \nonumber \\
&&\leq  \P_{\tau}[H^{-1}> \ell]+\E^\Q_{\tau}\left[H^{-1}\1_{\{H^{-1}\leq \ell\}}\1_{\{\tilde  \theta_{r}
       (\tau)<\theta_{0}(\tau)\leq \tau + T\}}\1_{\Bc_T}\right] \nonumber \\
&&\leq  \P_{\tau}[H^{-1}> \ell]+\ell\ \Q_{\tau}
       \left[\tilde  \theta_{r}(\tau)<\theta_{0}(\tau)\leq \tau + T, \Bc_T\right]. 
       \label{eq: decompo P Ac:3}
\end{eqnarray}
The first term above can be estimated, by using  Chebyshev's inequality, the inequality $\theta_{0,r}^T\le \theta_{0}(\tau)$, and Lemma \ref{lemma: Ecarre tau_new}:
\be
      \Pcond{\tau}{H^{-1}> \ell}
&\le& \frac{1}{|\log \ell|^{2}} {2 \Econd{\tau}{\frac{ 1 }{4 }
      \lambda^4_\infty|\theta_{0,r}^T-\tau|^{2} +\lambda_\infty^2(\theta_{0,r}^T-\tau)}}
      \nonumber \\
&\le& {\frac{(L\lambda^4_\infty+2\lambda^2_\infty)}{|\log \ell|^{2}}  \Phi^{1}_{0}(\tau)
\leq  \frac{ c }{2 } \Phi^{1}_{0}(\tau)},\label{eq: decompo P Ac:4_new}
\ee
where the last inequality holds by taking the constant $\ell=\ell(c,L,d)$ large enough.
To handle the second term in \reff{eq: decompo P Ac:3}, set 
$$
M_t:= \E^\Q_t[P_\tau]+\1_{[\tau,\infty)}(t)\int_\tau^t \1_{\{s<\theta_{0,r}^T\}} a_s^\top dW^\Q_s
\;\;\mbox{ for } \;\; t\ge 0  
$$
so that $M$ is a u.i. $\Q$-martingale.
Let $\theta^M_{0}(\tau)$ and $\theta^M_{r}(\tau)$ be the first hitting times after $\tau$ of levels 0 and $r$ by $M$, and set $\theta_{0,r}^{M,T}:=\theta^M_{0}(\tau) \wedge  \theta^M_{r}(\tau)\wedge (\tau+T).$ Recalling assumption \refAssOZx , we see that   
$$
  \1_{[\tau,\infty)}(t)\int_\tau^t \1_{\{s<\theta_{0,r}^T\}} a_s^\top \lambda_s ds
= \1_{[\tau,\infty)}(t)\int_\tau^t \1_{\{s<\theta_{0,r}^T\}} b_s  ds
  \;\;\mbox{ on } \Bc_T\cap E.
$$
Hence, on $\Bc_T\cap E$ the processes $M$ and $P$ coincide on 
$[\tau,\theta_{0,r}^T]$.
By the optional sampling theorem and the non-negativity of $(M_t)_{t\in [\tau,\theta_{0,r}^{M,T}]}$ 
a.s. on $E$, we then deduce
\b*
       P_\tau \1_E
  =    M_\tau \1_E
  =    \E^\Q_\tau(\1_E M_{\theta_{0,r}^{M,T}})
&\geq& \1_E r \Q_\tau\left(  \theta^M_{r}(\tau)<\theta^M_{0}(\tau)\land  (\tau + T)\right) \\
&\geq& \1_E r \Q_\tau\left(\tilde  \theta_{r}(\tau)<\theta_{0}(\tau)<  \tau + T, \ \Bc_T\right).
\e*
Plugging this inequality together with \reff{eq: decompo P Ac:4_new} 
into \reff{eq: decompo P Ac:3}, gives on $E$ that
$$
\P_{\tau}\left[\tilde  \theta_{r}(\tau)<\theta_{0}(\tau)\leq \tau + T, \Bc_T\right]   \leq  \frac{ c }{2 } \Phi^{1}_{0}(\tau)+\frac{\ell(c,L,d)}r P_\tau . \label{eq: decompo P Ac:6}
$$
\ep
\vs3

{\bf Proof of Theorem \ref{thm: bound on tauxa_new}.} 
\underline{Part (a):} For $\tau \in \Tc^{\pi}$ and $c>0$
Lemmas \ref{lem: first bound proof thm bound tauxa} and 
\ref{lem: bound P Axacomp} imply
$$
     \Phi^{1}_{0}(\tau) 
 \le c_{\reff{lem: first bound proof thm bound tauxa}}  {\Pcond{\tau}{\Ac^\tau_{0}}} \le c_{\reff{lem: first bound proof thm bound tauxa}} \left [ \eta(c) P_{\tau} + c\ \Phi^{1}_{0}(\tau) + \kappa\left (\frac{2}{c},r\right ) |\pi| \right ]
$$
a.s. on $\{P_\tau \in [0,r]\}$. {Specializing} to $c=1/(2c_{\reff{lem: first bound proof thm bound tauxa}})$ leads to
$$
      \Phi^{1}_{0}(\tau) 
 \le 2 c_{\reff{lem: first bound proof thm bound tauxa}} 
      \left [ \eta((2c_{\reff{lem: first bound proof thm bound tauxa}})^{-1}) P_{\tau}  
       + \kappa\left (4c_{\reff{lem: first bound proof thm bound tauxa}} ,r\right ) |\pi|\right ] $$
on $\{ P_\tau\in [0,r]\}$.
 On $\{ P_\tau <0 \}$ we simply have
$\Phi^{1}_{0}(\tau) = 0$, while  
$      \Phi_0^1(\tau)    
   \le L  $
on $\{ P_\tau > r \}$ by assumption \refassbornetaux .
{  
This implies 
\begin{equation}\label{eqn:thm: bound on tauxa_old}
 \Phi^{1}_{0}(\tau)
         \le   {   \bar c_{\reff{thm: bound on tauxa_new}} } \Big [ P_{\tau} + |\pi| \Big ] \1_{\{0\le  P_\tau\le r\}}
      +  L \1_{\{ r<P_\tau  \}},
\end{equation}
      where ${\bar c_{\reff{thm: bound on tauxa_new}} = \bar c_{\reff{thm: bound on tauxa_new}} }(r,L,d,\kappa)>0$. By a change of the constant ${\bar c_{\reff{thm: bound on tauxa_new}} }$ the assertion follows.
}
\vs2

\underline{Part (b):}
{
Combining \eqref{eqn:thm: bound on tauxa_old} and Lemma \ref{lemma:A_to_B}, we derive
$$ \Phi^{1}_{0}(\tau)
         \le \bar d_{\reff{thm: bound on tauxa_new}}  
         \Big [ |P_{\tau}| + |\pi|^\frac{1}{2} \Big ] 
             \1_{\{ |P_\tau|\le \tilde r\}}
             +  L \1_{\{ \tilde r<|P_\tau|  \}}$$ 
for $0<\tilde r<r$ and $\tau\in \Tc$, where
$  \bar d_{\reff{thm: bound on tauxa_new}} 
 = \bar d_{\reff{thm: bound on tauxa_new}}(r,\tilde r,L,d,\kappa)>0$. 
 {Observe that $\Phi^{1}_{0}(\tau)=0$ for $P_\tau\leq 0$, thus the above r.h.s. needs to be specialized only for $P_\tau\geq 0$. Then, }
 choosing $\tilde r = r/2$
and adapting $  \bar d_{\reff{thm: bound on tauxa_new}}$, we obtain part (b).
\ep}


\subsection{Proof of Theorem \ref{thm: bound on taupi0_new}}
\label{sec: proof of thm: bound on taupi0}
  
Now we are in a position to conclude the proof of  Theorem \ref{thm: bound on taupi0_new}.  It is based on a recursion argument. Namely, given $\tau\in \Tc$ such that $0\le P_{\tau}\le r$, we wait until the next time $\vartheta$ in $\R_{+}$ such that $Z$ hits the boundary. The time it takes, $\vartheta-\tau$, is controlled by  Theorem \ref{thm: bound on tauxa_new}. If $Z_{\phi^{+}_{\vartheta}}\notin \Oc$, then we stop: $\theta_{0}^{\pi}(\tau)-\tau\le \vartheta-\tau +|\pi|$. If not, then we know from standard estimates (Lemma \ref{lemma: ex ass iii} below) that $P_{\phi^{+}_{\vartheta}}\in [0,r]$,
up to some event with a probability controlled by $O(|\pi|^{\frac12})$. In this case 
one can restart the above procedure from $\phi^{+}_{\vartheta}\in \pi$. Again, one waits for the next time in $\R_{+}$ such that $Z$ reaches the boundary and stops if $Z\notin \Oc$ at the following time in $\pi$. One iterates this procedure.  The key point is that the probability of the event set $\{ Z_{\phi^{+}_{\vartheta}}\in \Oc\}$  is uniformly controlled by some $\alpha<1$ (see Lemma \ref{lem: proba cross on boundary} below). 
\vs3

Before we start with the proof of Theorem \ref{thm: bound on taupi0_new} we state two lemmas that are needed.
The first one can be verified by Doob's maximal inequality and 
assumption \refAssOZxi:
\vs3

\begin{lemma}\label{lemma: ex ass iii} 
Under the assumption  \refAssOZxi\  one has, for all $\tau \in \Tc$ and $\lambda>0$,
$$
    \Pcond{\tau}{\max_{\tau \le t\le \phi^{+}_{\tau}} |P_{t}-P_{\tau}|\ge \lambda}
\le \frac{1}{\lambda} \Econd{\tau}{\max_{\tau \le t\le \phi^{+}_{\tau}} |P_{t}-P_{\tau}|^{2}}^{\frac12}
\le   \frac{c_{(\ref{lemma: ex ass iii})}}{\lambda} |\pi|^{\frac12}, $$
where $c_{(\ref{lemma: ex ass iii})}:=L+{2}\sqrt{d}L$.
\end{lemma}
\vs3

\begin{lemma}\label{lem: proba cross on boundary} 
Let assumptions \refAssZx\  and \refAssOZx\ hold. Then there exists an  
$0<\alpha_{\reff{lem: proba cross on boundary} }=\alpha_{\reff{lem: proba cross on boundary} }(r,L,d,\kappa)<1$   such that, a.s.,
$$     \P_\tau\left[P_{\phi^{+}_{\tau}}>0\right]
   \le \alpha_{\reff{lem: proba cross on boundary} } 
   \;\; \mbox{ on } \;\;
   \left \{ \gamma(\tau,\phi_\tau) \le \frac{r}{2}, P_\tau=0 \right \} \in \Fc_\tau
   $$
for all $\tau\in \Tc$ and $0<|\pi|\le \eps_{\reff{lem: proba cross on boundary}} 
  = \eps_{\reff{lem: proba cross on boundary}} 
 (r,L,d,\kappa)$.

\end{lemma}
\vs3

\proof 
It is sufficient to check for $B\in \Fc_\tau$ of positive measure with 
$$B\subseteq \left \{ \gamma(\tau,\phi_\tau) \le \frac{r}{2}, P_\tau=0  \right \} $$
that 
$$     \Pro{P_{\phi^{+}_{\tau}}>0,B} 
   \le \alpha_{\reff{lem: proba cross on boundary} } \Pro{B}.$$
Let
$$    \Bc 
   := \left \{\max_{\tau \le t\le \phi^{+}_{\tau}} (|P_{t}|\vee \gamma({t},{\phi_{t}}\vee \tau))
            \le \frac{r}{2} \right \}
      \cap B $$
so that
\begin{equation}\label{eqn:B_in_good_set}
   \Bc \subseteq 
   \left \{\max_{\tau \le t\le \phi^{+}_{\tau}} (|P_{t}|\vee \gamma({t},{\phi_{t}}))
            \le r \right \}.
\end{equation}
We use assumptions \refAssOZx, \refAssZx\ and Lemma \ref{lemma: ex ass iii} 
to continue with
\b*
&   & \Pro{P_{\phi^{+}_{\tau}}>0, B} \\
& = & \Pro{\int_\tau^{\phi^+_{\tau}} b_{s} ds 
                 + \int_\tau^{\phi^+_{\tau}} a^{\top}_{s} dW_{s} > 0, B } \\
&\le& \Pro { \int_\tau^{\phi^+_{\tau}} a^{\top}_{s} dW_{s} > 
      - L(\phi^+_{\tau}-\tau), \Bc} + \Pro{\Bc^c \cap B} \\
&\le& \Pro{ \int_\tau^{\phi^+_{\tau}} a^{\top}_{s} dW_{s} > 
      - L(\phi^+_{\tau}-\tau), \Bc } \\  
&   & + \left [\kappa\left (1,\frac{r}{2}\right ) |\pi|+\frac{2c_{(\ref{lemma: ex ass iii})}}{r} |\pi|^{\frac12}\right ] \Pro{B}.
\e*
Assuming that we are able to show that
\begin{equation}\label{eqn:linear_diffusion_lower_bound}
 \Pro{ \int_\tau^{\phi^+_{\tau}} a^{\top}_{s} dW_{s} > 
      - L(\phi^+_{\tau}-\tau), \Bc }
   \le \theta \Pro{B}
\end{equation}
for some $\theta=\theta(L,d)\in (0,1)$, the proof would be complete as
$$     \Pro{P_{\phi^{+}_{\tau}}>0, B}
   \le \left [ \theta  + \kappa\left ( 1,\frac{r}{2}\right )|\pi|+
       \frac{2 c_{(\ref{lemma: ex ass iii})}}{r} |\pi|^{\frac12} 
       \right ] \Pro{B}$$
and $\eps_{\reff{lem: proba cross on boundary}} = \eps_{\reff{lem: proba cross on boundary}}(r,L,d,\kappa)>0$ can be taken small enough to {guarantee} 
$$     \Pro{P_{\phi^{+}_{\tau}}>0, B}
   \le \alpha \Pro{B}$$
for some $\alpha=\alpha(r,L,d,\kappa)\in (0,1)$.
In order to check \reff{eqn:linear_diffusion_lower_bound} we let
$$ M_t := e^{\top} W_{t \wedge \tau}+ \int_\tau^{\tau\vee t} \bar a_s^\top dW_s $$
where $\bar a_s := a_s \1_{\{s\le \phi_\tau^+ \}} + e \1_{\{s > \phi_\tau^+ \}}$
with $e=d^{-\frac12}(1,...,1)^{\top}$. 
Define $\Lambda(s):=\inf\{t\ge 0:  \langle M \rangle_{t}>s\}$. 
Applying the Dambis-Schwarz Theorem \cite[p. 181]{rev:yor:05} yields that
$B:= M_{\Lambda}$ is a Brownian motion in the filtration $\G=(\Gc_{t})_{t\ge 0}$ defined by $\Gc=\Fc_{\Lambda}$
and $M=B_{\langle M\rangle}$. One can also check that
$\Fc_\tau \subseteq \Gc_{\langle M \rangle_\tau}$.
Letting $\eta := \phi^+_{\tau}-\tau$ (which is $\Fc_\tau$-measurable), observing
$\eta\le \sqrt{\eta}$ and
$ \eta L^{-2} \le \langle M \rangle_{\tau+\eta}-\langle M \rangle_{\tau} \le \eta d  L^2 $ on $\Bc$ by assumption \refAssOZx\ and \reff{eqn:B_in_good_set}, and taking an auxiliary one-dimensional Brownian 
motion $\widetilde B$ defined on some $(\widetilde{\Omega},\widetilde{\P})$, we conclude by
\b*
&   & \Pro{     \int_\tau^{\phi^+_{\tau}} a^{\top}_{s} dW_{s} > 
                 - L(\phi^+_{\tau}-\tau), \Bc } \\
& = & \Pro{     B_{\langle M \rangle_{\phi_\tau^+}} 
                 - B_{\langle M \rangle_\tau} > 
                 - L(\phi^+_{\tau}-\tau), \Bc } \\
&\le& \Pro{    \sup_{t\in [\langle M \rangle_\tau +\eta L^{-2},
                               \langle M \rangle_\tau +\eta d L^2]}
                   B_t
                 - B_{\langle M \rangle_\tau} > 
                 - L\eta, \Bc } \\
&\le& \Pro{    \sup_{t\in [\langle M \rangle_\tau +\eta L^{-2},
                               \langle M \rangle_\tau +\eta d L^2]}
                   B_t
                 - B_{\langle M \rangle_\tau} > 
                 - L\eta, B  } \\
&\le& \widetilde{\P}\times \Pro{      \sup_{\eta L^{-2} \le u \le \eta  d L^2} \widetilde B_u 
                 >  -L \sqrt{\eta}, B } \\
&\le& \widetilde{\P}\times \Pro{      \sup_{L^{-2} \le u \le d L^2} \widetilde B_u 
                 >  -L , B } \\
& = & \widetilde{\P} \left [ \sup_{L^{-2} \le u \le d L^2} \widetilde B_u 
                 >  -L \right ] \Pro{B} \\
& =:& \theta \Pro{B}.
\e*
\ep

{\bf Proof of Theorem \ref{thm: bound on taupi0_new}.}  (a)
{
First we assume that $\tau\in \Tc^\pi$. For
}
$i\ge 0$ we define
$$ \vartheta_{0}   :=\theta_{0}(\tau) \;,\quad 
   \vartheta_{i+1} :=\theta_{0}(\phi^{ +}_{\vartheta_{i }})\;, \quad
   \vartheta^{\pi}_{0}:=\theta^{\pi}_{0}(\tau)\;,\quad 
   \vartheta^{\pi}_{i+1}:=\theta^{\pi}_{0}(\phi^{ +}_{\vartheta_{i }}), $$
$$
E_{i}:=\{P_{\phi^{+}_{\vartheta_{i}}}>0\},
\;\;\mbox{ and }\;\; \Ac_i:=\cap_{0\leq j\leq i }E_{j}
\in \Fc_{\phi^{+}_{\vartheta_{i}}}.
$$
1. From the definitions we obtain for $i\ge 0$:
\begin{enumerate}[a)]
\item $\vartheta_{i+1}\leq\vartheta^{\pi}_{i+1}$ (by definitions of the stopping times);
\item $\phi^{+}_{\vartheta_{i+1}}\leq\vartheta_{i+1}+|\pi|$ (by the definition of $\phi^{+}$);
\item   $\vartheta^{\pi}_{i+1}=\vartheta^{\pi}_{i+2}$ on $E_{i+1}=\{Z_{\phi^{+}_{\vartheta_{i+1}}}\in \Oc\}$   (since $\phi^{+}_{\vartheta_{i+1}}< \vartheta^{\pi}_{i+1}$ on $E_{i+1}$); 
\item   $\vartheta^{\pi}_{i+1}\leq \vartheta_{i+1}+|\pi|$  on $(E_{i+1})^c=\{Z_{\phi^{+}_{\vartheta_{i+1}}}\notin \Oc\}$ (by definition of the stopping time $\vartheta^{\pi}_{i+1}$).
\end{enumerate}
Item c) leads to
\begin{align*}
\vartheta^{\pi}_{i+1} &= \vartheta^{\pi}_{i+2}\1_{E_{i+1}}+ \vartheta^{\pi}_{i+1}\1_{[E_{i+1}]^c},
\end{align*}
$$
     \vartheta^{\pi}_{i+1}-\phi^{+}_{\vartheta_{i}}
 =   (\vartheta^{\pi}_{i+2}-\phi^{+}_{\vartheta_{i+1}})\1_{E_{i+1}}
    +(\phi^{+}_{\vartheta_{i+1}}-\phi^{+}_{\vartheta_{i}})\1_{E_{i+1}}
    +(\vartheta^{\pi}_{i+1}-\phi^{+}_{\vartheta_{i}})\1_{[E_{i+1}]^c}.
$$
With b) and d) we continue to
\b*
\vartheta^{\pi}_{i+1}-\phi^{+}_{\vartheta_{i}}
&\leq& (\vartheta^{\pi}_{i+2}-\phi^{+}_{\vartheta_{i+1}})\1_{E_{i+1}}+(
\vartheta_{i+1}+|\pi|-\phi^{+}_{\vartheta_{i}})\1_{E_{i+1}}\\
&& + ( \vartheta_{i+1}+|\pi|-\phi^{+}_{\vartheta_{i}}  )\1_{[E_{i+1}]^c} \\
& = & (\vartheta^{\pi}_{i+2}-\phi^{+}_{\vartheta_{i+1}})\1_{E_{i+1}}+|\pi|+(
\vartheta_{i+1}-\phi^{+}_{\vartheta_{i}})
\e*
and 
\b*
\Econd{\tau}{(\vartheta^{\pi}_{i+1}-\phi^{+}_{\vartheta_{i}})\1_{ \Ac_i }}
&\leq& \Econd{\tau}{(\vartheta^{\pi}_{i+2}-\phi^{+}_{\vartheta_{i+1}})\1_{ \Ac_i }\1_{E_{i+1}}}
\\
&&+|\pi|\Pcond{\tau}{\Ac_i}+
\Econd{\tau}{(\vartheta_{i+1}-\phi^{+}_{\vartheta_{i}})\1_{ \Ac_i }}\\
&=& \Econd{\tau}{(\vartheta^{\pi}_{i+2}-\phi^{+}_{\vartheta_{i+1}})\1_{ \Ac_{i+1} }}
\\
&&+|\pi|\Pcond{\tau}{\Ac_i}+
\Econd{\tau}{(\vartheta_{i+1}-\phi^{+}_{\vartheta_{i}})\1_{ \Ac_i }}.
\e*
Summing up the above inequalities from $i=0$ to $i=n-1$ yields
\be
\Econd{\tau}{(\vartheta^{\pi}_{1}-\phi^{+}_{\vartheta_{0}})\1_{ \Ac_0 }}
&\leq& \Econd{\tau}{(\vartheta^{\pi}_{n+1}-\phi^{+}_{\vartheta_{n}})\1_{ \Ac_n }}
\nonumber\\
&&+\sum_{i=0}^{n-1}\left(|\pi|\Pcond{\tau}{\Ac_i}+
\Econd{\tau}{(\vartheta_{i+1}-\phi^{+}_{\vartheta_{i}})\1_{ \Ac_i }}\right).\label{eq: control after summing up}
\ee
2. For $\sigma \in \Tc$ set 
$\Ac^\sigma:=\{\gamma(\sigma,\phi_{\sigma})\leq r/2 \}\in \Fc_\sigma$ so that, 
for $i\geq 1$,
\b*
      \Pcond{\tau}{\Ac_i}
& = & \Econd{\tau}{\1_{\Ac_{i-1}} \1_{\Ac^{\vartheta_{i}}}\Pcond{\vartheta_{i}}{E_{i}}} +  
      \Econd{\tau}{\1_{\Ac_{i-1}} \Pcond{\phi^{+}_{\vartheta_{i-1}}}
      {[\Ac^{\vartheta_{i}}]^c\cap E_{i}}} \\
&\le& \alpha_{\reff{lem: proba cross on boundary} } \Pcond{\tau}{\Ac_{i-1}}+
 \Econd{\tau}{\1_{\Ac_{i-1}} \Pcond{\phi^{+}_{\vartheta_{i-1}}}
      {[\Ac^{\vartheta_{i}}]^c\cap E_{i}}} ,
\e*
because of $\Fc_{\phi_{\vartheta_{i-1}}^+}\subseteq \Fc_{\vartheta_i}$, Lemma \ref{lem: proba cross on boundary},  and $P_{\vartheta_i}=0$ on $\Ac_{i-1}$.
To treat the second term we take a fixed $T>0$ and use \refAssZx\ and \refassbornetaux\ to get
\b*
      \Pcond{\phi^{+}_{\vartheta_{i-1}}}{[\Ac^{\vartheta_{i}}]^c\cap E_{i}}
&\le& \Pcond{\phi^{+}_{\vartheta_{i-1}}}{\{\gamma(\vartheta_{i},\phi_{\vartheta_{i}}) > r/2 \}
      \cap \{\vartheta_{i}\le \phi^{+}_{\vartheta_{i-1}}+T\}} \\
&   & + \Pcond{\phi^{+}_{\vartheta_{i-1}}}{ \vartheta_{i} > \phi^{+}_{\vartheta_{i-1}}+T}\\
&\le& \Pcond{\phi^{+}_{\vartheta_{i-1}}}{ \sup_{\phi^{+}_{\vartheta_{i-1}}\le t \le 
      \phi^{+}_{\vartheta_{i-1}}+T }\gamma(t,\phi_{t}\vee\phi^{+}_{\vartheta_{i-1}} ) > r/2 } \\
&   & + \Pcond{\phi^{+}_{\vartheta_{i-1}}}{ \vartheta_{i} > \phi^{+}_{\vartheta_{i-1}}+T}\\
&\le& \kappa(T,r/2)|\pi| + \E_{\phi^{+}_{\vartheta_{i-1}}}[\theta_{0}(\phi^{+}_{\vartheta_{i-1}})]/T\\
&\le& \kappa(T,r/2)|\pi| + L/T.
\e*
By taking $T>0$ large enough and then 
$   \varepsilon_{\reff{eq: control PAci by alpha}}\in 
    (0,\varepsilon_{\reff{lem: proba cross on boundary}}]$ small enough such that
$$  \alpha_{\reff{lem: proba cross on boundary} }
  + \kappa(T,r/2)\varepsilon_{\reff{eq: control PAci by alpha}}
  + \frac{L}{T} =:\alpha<1$$ 
 and assuming that 
$|\pi|\leq \varepsilon_{\reff{eq: control PAci by alpha}}$, we obtain
$\Pcond{\tau}{\Ac_i} \le \alpha \Pcond{\tau}{\Ac_{i-1}}$ and, by induction,
\begin{equation}\label{eq: control PAci by alpha}
       \Pcond{\tau}{\Ac_j}
\leq \alpha^{{j}}\;\;\pourtout\; j\ge 0.
\end{equation}
3. Let us set $F_{i}:=\{P_{\phi^{+}_{\vartheta_{i}}} > r\}$ for $i\ge 0$. Because of
$\phi^{+}_{\vartheta_{i}}\in \Tc^\pi$, applying
{ 
\eqref{eqn:thm: bound on tauxa_old} from the proof of Theorem \ref{thm: bound on tauxa_new}
}
 and using the fact that $\Ac_{i}\in \Fc_{\phi^{+}_{\vartheta_{i}}}$ and assumption \refassbornetaux, lead to 
\b*
&    & \Econd{\tau}{(\vartheta_{i+1}-\phi^{+}_{\vartheta_{i}})\1_{ \Ac_i }} \\
& =  & \Econd{\tau}{\Econd{\phi^{+}_{\vartheta_{i}}}{(\vartheta_{i+1}
       - \phi^{+}_{\vartheta_{i}})\1_{ \Ac_{i} }}} \\
&\leq& \Econd{\tau}{
       \bar c_{\reff{thm: bound on tauxa_new}} 
       \big(
       (P_{\phi^{+}_{\vartheta_{i}}})^+ + |\pi|\big)\1_{ \Ac_{i-1} }\1_{E_{i}\cap 
       [F_{i}]^{c}}} + L \Pcond{\tau}{\Ac_{i} \cap  F_{i} }\\
&\leq& \Econd{\tau}{
       \bar c_{\reff{thm: bound on tauxa_new}} 
       \big(
       (P_{\phi^{+}_{\vartheta_{i}}})^++|\pi|\big)\1_{ \Ac_{i-1} } }
+ L {\Econd{\tau}{1_{\Ac_{i-1}}\Pcond{\vartheta_{i}}{ 
  F_{i}} }},
\e*
where $\Ac_{-1}:=\Omega$.
Because $P_{\vartheta_{i}}\le 0$,  Lemma \ref{lemma: ex ass iii} implies 
$$     \Econd{\vartheta_i}{(P_{\phi_{\vartheta_i}^+})^+}
   \le c_{\reff{lemma: ex ass iii}} |\pi|^\frac{1}{2}
   \;\;\mbox{ and } \;\;
   \Pcond{\vartheta_i}{F_i} \le \frac{c_{\reff{lemma: ex ass iii}}}{r} |\pi|^\frac{1}{2}, $$
and \reff{eq: control PAci by alpha} yields
\b*
&    & \Econd{\tau}{(\vartheta_{i+1}-\phi^{+}_{\vartheta_{i}})\1_{ \Ac_i }} \\
&\leq& \Econd{\tau}{
       \bar c_{\reff{thm: bound on tauxa_new}} \big( c_{\reff{lemma: ex ass iii}}
       |\pi|^{\frac{1}{2}}+|\pi|^{\frac12}\big)\1_{ \Ac_{i-1} }}
       + L \Pcond{\tau}{\Ac_{i-1}} \frac{c_{(\ref{lemma: ex ass iii})}}{r} |\pi|^{\frac12} \\
&\leq& D |\pi|^{\frac{1}{2}} \alpha^{(i-1)_+}
\e*
with $D:=  \bar c_{\reff{thm: bound on tauxa_new}} c_{\reff{lemma: ex ass iii}}
         + \bar c_{\reff{thm: bound on tauxa_new}} 
         + L c_{(\ref{lemma: ex ass iii})}/r$.
If we insert the last estimate into \reff{eq: control after summing up} and let $n\to +\infty$,
then we get
$$
     \Econd{\tau}{{(\vartheta^{\pi}_{1}-\phi^{+}_{\vartheta_{0}})\1_{ \Ac_0 }}}
\leq |\pi|^{\frac{1}{2}}\frac{ |\pi|^{\frac{1}{2}}+(2-\alpha) D}{1-\alpha },
$$
where we exploit Lemma \ref{lemma: Ecarre tau_new} to check
$$ \Econd{\tau}{|\vartheta_{n+1}^\pi-\phi_{\vartheta_n}^+| \1_{\Ac_n}}
   \le \sqrt{L^{(2)}} \Pcond{\tau}{\Ac_n}^\frac{1}{2}.$$
Observe now that 
$$    \theta^{\pi}_{0}(\tau)
    = \big[ \phi^{+}_{\vartheta_{0}}+(\vartheta^{\pi}_{1}-\phi^{+}_{\vartheta_{0}})\big]
      \1_{ \Ac_0 } +\phi^{+}_{\vartheta_{0}}\1_{ [\Ac_0]^c }
  \leq |\pi|+\vartheta_{0}+(\vartheta^{\pi}_{1}-\phi^{+}_{\vartheta_{0}})\1_{ \Ac_0 },$$
so that by an application of the previous estimate, 
 \eqref{eqn:thm: bound on tauxa_old} and Assumption {\bf (L)} we obtain  
\begin{equation}\label{eqn:thm: bound on taupi0_old}
\Phi^{1,\pi}_{0}(\tau) \le    \bar c_{\reff{thm: bound on taupi0_new} }\left[ P_{\tau} +|\pi|^{\frac12}  \right]\1_{\{0\le  P_\tau\le   r\}}
             +  L \1_{\{   r<P_\tau  \}},
\end{equation}
for $\tau\in \Tc^{\pi}$ and 
$   \bar c_{\reff{thm: bound on taupi0_new} }
  = \bar c_{\reff{thm: bound on taupi0_new} }(r,L,d,\kappa)>0$.
  
(b) We now consider the general case $\tau \in \Tc$.
Applying Lemma \ref{lemma:A_to_B} to \eqref{eqn:thm: bound on taupi0_old}
we obtain for $0<\tilde r<r$ that
{\b*
\Phi^{1,\pi}_{0}(\tau) \le    \bar d_{\reff{thm: bound on taupi0_new} }\left[ |P_{\tau}| +|\pi|^{\frac12}  \right] \1_{\{|P_\tau|\le \tilde r\}}
             +  L \1_{\{ \tilde r<|P_\tau|  \}},
\e* }
where 
$   \bar d_{\reff{thm: bound on taupi0_new} }
  = \bar d_{\reff{thm: bound on taupi0_new} }(r,\tilde r,L,d,\kappa)>0$.
Taking $\tilde r = r/2$ and adapting $\bar d_{\reff{thm: bound on taupi0_new} }$,
we obtain the statement of the theorem.
\ep


\section{General ${\rm \bf L}_{1}$-error for exit time approximations}
\label{section:L1}

The main application we develop in this paper is the study of the error made by estimating the exit time $\theta$ of a diffusion $X$ from a domain $\Oc$ by the discrete exit time $\bar \theta$ of  an approximation process $\bar X$, which can be $X$ itself or its Euler or Milstein scheme etc, computed on a grid $\bar \pi$. We only assume that the corresponding distance processes remain close, at least at the order $|\bar \pi|^{\frac12}$ in ${\rm \bf L}_{1}$. 
 If $X$ exits before $\bar X$, then 
 our
assumptions imply that $\bar X$ is close to the boundary as well. If we also know that the expectation of the time it takes to the approximation scheme $\bar X$ to exit the domain is proportional to its distance to the boundary up to an additional term $|\bar \pi|^{\frac12}$, then we can conclude that $\E[|\bar \theta-\theta|\1_{\{\theta\le \bar \theta\}}]$ is controlled in $|\bar \pi|^{\frac12}$. The same idea applies if $\bar X$ exits before $X$. 
In this section, we show how Theorems \ref{thm: bound on tauxa_new} and \ref{thm: bound on taupi0_new} 
are used to follow this idea. We start with an abstract statement and then specialize it to the case where $X$ solves a stochastic differential equation and $\bar X$ is its Euler scheme.


\subsection{Upper-bound in an abstract setting} \label{sec: upper bound abstract setting}
We fix an open non-empty subset $\Oc$ of a metric space $(\Zc,d_{\Zc})$, satisfying the assumptions of 
Section \ref{subsec: pbm formulation abstract}, and two {$\Zc$-valued} processes $X$ and $\bar X$.
We consider the first exit time $\theta_{0}:=\theta_{0}(0)$ of $X$ on $\pi:=\R_{+}$ and $\bar \theta^{\bar \pi}_{0}:=\bar \theta^{\bar \pi}_{0}(0)$  of $\bar X$ on $\bar \pi\subsetneq \R_{+}$  
(where $\bar \pi$ satisfies the conditions of 
Section \ref{subsec: pbm formulation abstract}), i.e.
$$
     \theta_{0}
  := \inf\{t\ge 0~:~X_{t}\notin \Oc\}
     \;\; \mbox{ and }\;\; 
     \bar \theta^{\bar \pi}_{0}
     :=\inf\{t\ge 0~:~t\in \bar \pi \mbox{ and } \bar X_{t}\notin \Oc\}.$$
We let $\bar \phi$ and $\bar \phi^{+}$ be {the functions defined in \eqref{eq:phi:phiplus}} associated to $\bar \pi$.
\vs2

We also fix a     {\sl distance function}  $\delta:\Zc\mapsto \R$ such that $\delta>0$ on $\Oc$, 
$\delta=0$   on $\partial \Oc$, and  $\delta<0$ on $\bar \Oc^{c}$, and  set $P:=\delta(X)$ and $\bar P:=\delta(\bar X)$.  
\vs2

Throughout this section we assume that the assumptions  \refAssZx,   \refAssOZx\  and \refassbornetaux\ of Section \ref{subsec: pbm formulation abstract} hold for  $(X,\pi,P)$, 
$(\bar X, \bar \pi,\bar P)$, and $\delta$ with the same $(r,L,\kappa)$. Obviously, the estimate contained in  \refAssZx\ is trivial for $(X,\pi,P)$ since $\pi=\R_{+}$ and $\phi$ is the identity. 
 
\vs3

\begin{theorem}\label{thm: error exit times abstract_new}  Assume a stopping time $\upsilon:\Omega\to [0,\infty]$ and some $\rho>0$ such that
\be\label{eq: ass error stop time abstract}
    \Esp{|  P_{\vartheta}-\bar P_{\vartheta} |}
\le { \rho|\bar \pi|^{\frac12}}\;\;\;\mbox{ for all } \;\vartheta \in \Tc 
    \; \mbox{ with } \;
    \vartheta\le  \theta_0 \wedge \upsilon.
\ee 
Then for all integers $p\ge 1$ there exist  
$   c_{\reff{thm: error exit times abstract_new}}
  = c_{\reff{thm: error exit times abstract_new}}(r,L,d,\kappa,p,\rho)>0$ and 
$   \eps_{\reff{thm: error exit times abstract_new}}
  = \eps_{\reff{thm: error exit times abstract_new}}(r,L,d,\kappa)>0$ such that,
for $|\bar \pi|\le \eps_{\reff{thm: error exit times abstract_new}}$,
\[
        \Esp{  |[\theta_{0}\wedge \upsilon] -[\bar \theta^{\bar \pi}_{0}\wedge \upsilon]|^p}
    \le c_{\reff{thm: error exit times abstract_new}}\;   |\bar \pi|^{\frac12}.
\]
\end{theorem}

\proof Define
$\upsilon_0 :=  \theta_{0} \wedge \upsilon$
and
$\bar\upsilon_0 := \bar \theta^{\bar \pi}_{0} \wedge \upsilon$.
We observe that
\b*
      \Econd{\bar\upsilon_0}{   [\upsilon_{0}-\bar\upsilon_{0}]^p }
&\le& (\Phi_0^p(\bar\upsilon_0))^p \; \;\;\;\;
      \mbox{ on } \{\upsilon_{0}\ge  \bar\upsilon_{0}\},\\
      \Econd{\upsilon_0}{   [\bar\upsilon_{0}-\upsilon_{0}]^p}
&\le& (\Phi_0^{p,\pi}(\upsilon_0))^p \;\;\mbox{ on } 
      \{\upsilon_{0}<  \bar\upsilon_{0}\}
\e*
and continue with Lemma \ref{lemma: Ecarre tau_new} to get
\b*
      \Econd{\bar\upsilon_0}{   [\upsilon_{0}-\bar\upsilon_{0}]^p }
&\le& p! L^{p-1} \Phi_0^1(\bar\upsilon_0) \; \;\;\;\;
      \mbox{ on } \{\upsilon_{0}\ge  \bar\upsilon_{0}\},\\
      \Econd{\upsilon_0}{   [\bar\upsilon_{0}-\upsilon_{0}]^p}
&\le& p! L^{p-1} \Phi_0^{1,\pi}(\upsilon_0) \;\;\mbox{ on } 
      \{\upsilon_{0}<  \bar\upsilon_{0}\}.
\e*
Applying Theorem \ref{thm: bound on tauxa_new} to $(X,\pi,P)$ and
$\tau=\bar\upsilon_{0}$ we get
\b* 
        \Econd{  {\bar\upsilon_{0}}      }{  [\upsilon_{0}-\bar\upsilon_{0}]^p }
   &\le&  p! L^{p-1} c_{\reff{thm: bound on tauxa_new}}  P_{\bar\upsilon_{0}}
          \1_{\{P_{\bar\upsilon_{0}}\ge 0\}} \\
   &\le&  p! L^{p-1} c_{\reff{thm: bound on tauxa_new}} 
         |P_{\bar\upsilon_{0}}- \bar P_{\bar\upsilon_{0}}|
\e*
on $\{\upsilon_{0}> \bar\upsilon_{0}\}$, where we use that on  
$\{\upsilon_{0} >  \bar\upsilon_{0}\}$ we have
$\bar\upsilon_{0}=\bar \theta^{\bar \pi}_{0}$ and therefore 
$\bar P_{\bar\upsilon_{0}} = \bar P_{\bar \theta^{\bar \pi}_{0}}\le 0$. 
Consequently, 
\[     \Econd{  {\bar\upsilon_{0}}      }{  [\upsilon_{0}-\bar\upsilon_{0}]^p }
   \le  p! L^{p-1} c_{\reff{thm: bound on tauxa_new}} 
         |P_{\bar\upsilon_{0}}- \bar P_{\bar\upsilon_{0}}| 
   \mbox{ on }
   \{\upsilon_{0} \ge \bar\upsilon_{0}\}. \]
Applying 
Theorem \ref{thm: bound on taupi0_new} to $(\bar X,\bar \pi,\bar P)$ and $\tau=\upsilon_{0}$ implies
\b*
        \Econd{\upsilon_0}{  [ \bar\upsilon_{0}-\upsilon_{0}]^p}
    &\le& {p! L^{p-1}  d_{\reff{thm: bound on taupi0_new} } 
        [|\bar P_{\upsilon_{0}}| + |\bar \pi|^{\frac12}] }\\
&=&  p! L^{p-1} d_{\reff{thm: bound on taupi0_new} } 
        [|\bar P_{\upsilon_{0}} - P_{\upsilon_0}| + |\bar \pi|^{\frac12}]
\e*
on $\{\upsilon_{0}<  \bar\upsilon_{0}\}$, where (similarly as above) on this set 
$\upsilon_{0}=\theta_{0}$ and therefore 
$P_{\upsilon_{0}} = P_{\theta_{0}} = 0$. Letting $\vartheta := \upsilon_0\wedge  \bar\upsilon_{0}$,
the above inequalities imply
\[
      \Econd{{\vartheta}}{ |  \upsilon_{0}-\bar\upsilon_{0} |^p} 
 \le p! L^{p-1} [ c_{\reff{thm: bound on tauxa_new} } \vee
                  d_{\reff{thm: bound on taupi0_new} } ]
       \Big [   |   P_{\vartheta}- \bar P_{\vartheta} |
         + |\bar \pi|^{\frac12} \Big ],
\]
which, by Assumption \reff{eq: ass error stop time abstract}, leads to the desired result.
\ep

\medskip

We conclude this section with  sufficient conditions ensuring that \reff{eq: ass error stop time abstract} holds. 
In the following, $\|\cdot\|_{q}$ denotes the ${\rm \bf L}_{q}$-norm for $q\ge 1$. The proof being standard, it is 
postponed to the Appendix.

\begin{lemma}
\label{lem: error euler scheme up to stopping time with expo decay_new_new}
Assume $\vartheta,\tau\in \Tc$ such that $0\le \vartheta \le \tau$.
\begin{enumerate}[{\rm (a)}]
\item We have that
$$    \left \|\Dist{X_{\vartheta},\bar X_{\vartheta}} \right \|_1
  \le \inf_{1<q<\infty}   \sum_{k=0}^\infty \Pro{\tau \ge k}^\frac{q-1}{q}
                                \left \| \sup_{t\in [k,k+1)} \Dist{X_t,\bar X_t} \right \|_q
                        . $$
\item 
      Assume $\alpha>0$, $0 <\beta<\infty$, $1<q<\infty$, and
      $Q(\cdot,q):\{ 0,1,2,... \} \to \R_{+}$ such that
      \begin{enumerate}[{\rm (i)}]
      \item $\Pro{\tau\ge k} \le \alpha e^{-\beta k}$  for $k=0,1,2,...$,
      \item $\left \| \sup_{t\in [k,k+1)} \Dist{X_t,\bar X_t} \right \|_q
             \le Q(k,q) |\bar \pi|^\frac{1}{2}$ for $k=0,1,2,...$
      \item $ c:= \sum_{k=0}^\infty e^{\beta \left ( \frac{1}{q} -1 \right ) k}
                                 Q(k,q)<\infty.$
      \end{enumerate}
       Then one has 
      $      \left \|\Dist{X_{\vartheta},\bar X_{\vartheta}} \right \|_1
         \le \alpha^{1 - \frac{1}{q}} c |\bar\pi|^\frac{1}{2}=O(|\bar\pi|^\frac{1}{2})$.
\end{enumerate}
\end{lemma}
\medskip

Here we have some kind of trade-off between the decay of $\Pro{\tau\ge k}$, measured by $\beta$, and the 
growth of $\left \| \sup_{t\in [k,k+1)} \Dist{X_t,\bar X_t} \right \|_q$ measured by $Q(\cdot,\cdot)$.
In the product $e^{\beta \left ( \frac{1}{q} -1 \right ) k} Q(k,q)$ the factor 
$Q(k,q)$ is thought to be increasing in $q$, but the factor 
$e^{\beta \left ( \frac{1}{q} - 1\right ) k}$ decreases as $\beta$ and $q$ increase.

\vs2

Combining Theorem \ref{thm: error exit times abstract_new} and 
Lemma \ref{lem: error euler scheme up to stopping time with expo decay_new_new},
and using 
\[     |P_\theta-\bar P_\theta| = |\delta(X_\theta) - \delta(\bar X_\theta)|
   \le L \Dist{X_\theta,\bar X_\theta}, \]
gives  the following corollary. 

\begin{corollary}\label{cor: error exit times abstract_new}
Let $\upsilon$ be a stopping time. Assume that the conditions of Lemma \ref{lem: error euler scheme up to stopping time with expo decay_new_new}(b) are satisfied with 
$\tau=\theta_{0}{  \wedge\upsilon}$,
and let $p\ge 1$ be an integer. Then there 
is a $c>0$, depending at most on $(r,L,d,\kappa,p,{  \alpha,\beta,q,Q})$, such that 
\b*
    \Esp{  |[\theta_{0}{  \wedge \upsilon}] - [\bar \theta^{\bar \pi}_{0}{  \wedge \upsilon}]|^p} 
\le c\;   |\bar \pi|^{\frac12}\;\;\;\mbox{ whenever } \;\; |\bar \pi|\le   
    \eps_{\reff{thm: error exit times abstract_new}}
\e*
with $ \eps_{\reff{thm: error exit times abstract_new}}>0$ taken from 
Theorem \ref{thm: error exit times abstract_new}.
\end{corollary}
 

\subsection{Application to the Euler scheme approximation of the first exit time of a SDE}\label{sec: application to SDE}

Now we specialize the discussion to the case where $\Zc=\R^{d}$ endowed with the usual Euclidean norm $|\cdot|$ and where $X$ is the strong solution of the stochastic differential equation
\b*
X_{t}=x_{0}+\int_{0}^{t} \mu(X_{s}) ds + \int_{0}^{t} \sigma(X_{s}) dW_{s}
\e*
for some fixed $x_{0}\in \Oc$, where $(\mu,\sigma):\R^{d}\to (\R^{d},\R^{d\times d})$
satisfy 
\begin{assumption}\label{ass: coeff mu sigma}
There exists $0<L_\mu,L_\sigma \le L$ such that, for all $x,y \in \R^{d}$,
$$ |\mu(x)-\mu(y)|       \le L_\mu |   x-y|, \;\;\;\;
   |\sigma(x)-\sigma(y)| \le L_\sigma |x-y|, $$
and $|\mu(x)|+|\sigma(x)|\le L$.
\end{assumption}

\begin{remark}\label{rem : coeff sde} As usual some rows or columns  of $\sigma$ can be equal to $0$. In particular, the first component of $X$ can be seen as  the time component by setting the first entry of $\mu$ equal to $1$ and the first row of $\sigma$ equal to $0$, i.e. $X_t=(t,X^\flat_t)$ where $X^\flat$ is a diffusion process in $\R^{d-1}$. This allows to consider time dependent coefficients
(where one could investigate to what {extent} weaker assumptions on the first 
coordinate of $\mu$ and $\sigma$, like $1/2$-H\"older continuity, would be sufficient for the purpose of this paper).
This formalism allows also to consider time-dependent domains as in \cite{gobe:meno:10}, i.e. $\Oc=\bigcup_{t\geq 0} \left(\{t\}\times \Oc^\flat_t\right)$ where $(\Oc^\flat_t)_{t\geq 0}$ is a family of domains in $\R^{d-1}$. Then the distance function $\delta((t,x^\flat))$ shall be the signed spatial distance to the boundary $\Oc^\flat_t$. 
\end{remark}

In the following we denote by $D\delta$ 
and $D^{2}\delta$ the gradient (considered as row vector) and the Hessian matrix of 
$\delta$, respectively.  To verify condition \refAssOZx\, we use the following 
sufficient assumption:
\vs2

\begin{assumption}\label{ass: non characteristic boundary condition} 
There exists a bounded $C^{2}_{b}$ function $\delta:\R^d\mapsto \R$ such that 
$\delta>0$ on $\Oc$, $\delta=0$   on $\partial \Oc$ and  $\delta<0$ on $\bar \Oc^{c}$, which satisfies $|D\delta|\le 1$ and the non-characteristic boundary condition   
\be\label{eq: ass non chara boundary}
|D\delta\;\sigma| \ge {2}L^{-1}  \;\mbox{ on } \{|\delta|\le r\}.
\ee
\end{assumption}
 
Note that this condition is usually satisfied if $\sigma$ is uniformly elliptic and the domain  has a $C^{2}$ compact boundary, see e.g.  \cite{GiTr01}.
 \vs2

We let $\bar X$ be the Euler scheme based on the grid $\bar \pi$, i.e. 
\be\label{eq: def schema euler}
\bar X_{t}=x_{0}+\int_{0}^{t} \mu(\bar X_{\bar \phi_{s}}) ds + \int_{0}^{t} \sigma(\bar X_{\bar \phi_{s}}) dW_{s}.
\ee

We are now in a position to state the main results of this section, whose proofs are postponed to the end of the section. Note that a sufficient condition for the assumption \reff{eq: condition integrability for thm summary for sde} below is given in Lemma \ref{lemma:both_conditions_L}. See also \cite[Chapter 3]{frei:85}.
\medskip

\begin{theorem}\label{thm:summary_for_sde} 
Let the Assumptions \ref{ass: coeff mu sigma} and \ref{ass: non characteristic boundary condition} hold and assume that 
\be\label{eq: condition integrability for thm summary for sde}
 \Econd{\tau}{|\bar \theta^{\bar\pi}_{0}(\tau)-\tau|+|\theta_{0}(\tau)-\tau|} \le L 
   \;\;\mbox{ for all } \;\; \tau\in \Tc. 
   \ee
Let $\upsilon$ be a stopping time with values in $\R_{+}\cup \{\infty\}$.
   Assume that there are $\rho>0$, $4\le q <\infty$, and 
$\beta > \frac{qd}{q-1}(6 L_\mu+3q L_\sigma^2)$ such that 
$$ \Pro{\theta_{0}{  \wedge \upsilon} \ge k} \le \rho e^{-\beta k}
   \;\;\mbox{ for all } k=0,1,2,...$$
Then there exist $c,\eps>0$ and, for any integer $p\ge 1$, a constant $c_p>0$
such that, for  $|\bar \pi|\le \eps$,
$$ \Esp{  |[\theta_{0}{  \wedge \upsilon}]-[\bar \theta^{\bar \pi}_{0}{  \wedge \upsilon}]|^p}\le  c_p \;   |\bar \pi|^{\frac12}\;\;        
   \;\;\mbox{ and } \;\;
   \left(\Esp{  |X_{\theta_{0}{  \wedge \upsilon} }-\bar X_{ {\bar \theta^{\bar \pi}_{0}}{  \wedge \upsilon}}|^2}\right)^{\frac{ 1 }{2 }}
   \le c \;   |\bar \pi|^{\frac14}. $$
\end{theorem}
\medskip

\begin{remark}\label{rem:conditions_thm:summary_for_sde_new unbounded horizon}
Assuming (for example) $\upsilon\equiv \infty$, for the purpose of this paper the estimate 
$\Esp{  |\theta_{0} -\bar \theta^{\bar \pi}_{0}|^p}\le  c_p |\bar \pi|^{\frac12}$ is 
sufficient, as we know from \cite{gobe:00,gobe:meno:10} that it can not be improved 
for $p=1$. 
      However, it would be of interest to find the optimal exponents $\alpha_p>0$ such that 
      $\Esp{  |\theta_{0}-\bar \theta^{\bar \pi}_{0}|^p}\le  c_p |\bar \pi|^{\alpha_p}$, in the case $p>1$. This is left for future studies.
 \end{remark}

 In the case where we are only interested in a finite horizon problem, then the integrability condition 
 \reff{eq: condition integrability for thm summary for sde} is not necessary. 
 \medskip
  
\begin{theorem}\label{thm:summary_for_sde T finite} 
Let the Assumptions \ref{ass: coeff mu sigma} and \ref{ass: non characteristic boundary condition} hold.    
Fix $T>0$. Then there exist $c,\eps>0$ and, for any integer $p\ge 1$, a constant $c_p>0$
such that, for  $|\bar \pi|\le \eps$,
$$ \Esp{  |[\theta_{0}{  \wedge T}]-[\bar \theta^{\bar \pi}_{0}{  \wedge T}]|^p}\le  c_p \;   |\bar \pi|^{\frac12}\;\;        
   \;\mbox{ and } \;
   \left(\Esp{  |X_{\theta_{0}{  \wedge T} }-\bar X_{ {\bar \theta^{\bar \pi}_{0}}{  \wedge T}}|^2}\right)^{\frac{ 1 }{2 }}
   \le c \;   |\bar \pi|^{\frac14}. $$
\end{theorem}
\medskip

\begin{remark}\label{rem:conditions_thm:summary_for_sde_new finite T}
 The main aim of \cite{bouc:meno:09} was to study the strong  error made when approximating the solution of a BSDE whose 
      terminal condition is of the form $g(X_{\theta_{0}\wedge T})$, for some Lipschitz map $g$ and $T>0$, by a backward Euler 
      scheme; see \cite{bouc:meno:09} for the corresponding definitions and references. 
           Theorem \ref{thm:summary_for_sde T finite}  complements \cite[Theorem 3.1]{bouc:meno:09} in which the upper-bound takes the form 
      $O_{{|\bar \pi|\to 0}}(|\bar \pi|^{\frac12 - \eps})$ for all $0<\eps< 1/2$. Moreover, the upper-bound of the second inequality 
      of \cite[Theorem 3.3]{bouc:meno:09} is of the form $O_{{|\bar \pi|\to 0}}(|\bar \pi|^{\frac14 - \eps})$ for all $0<\eps<1/4$. 
      This comes from the control they obtained on the exit time of their Theorem 3.1. With Theorem \ref{thm:summary_for_sde T finite}  of this 
      paper it can be reduced to $O_{{|\bar \pi|\to 0}}(|\bar \pi|^{\frac14})$.  Our results open the door to the study of backward 
      Euler type approximations of BSDEs with a terminal condition  of the form $g(X_{\theta_{0}})$, i.e. there is no finite time 
      horizon $T>0$. This will however require to study at first the regularity of the solution of the BSDE, which is beyond  the scope 
      of this paper.
\end{remark}


\noindent {\bf Proof of Theorems \ref{thm:summary_for_sde} and  \ref{thm:summary_for_sde T finite}.}
(a) Theorem \ref{thm:summary_for_sde} is an immediate consequence of  Lemmas 
\ref{lem: error euler scheme up to stopping time with expo decay_new_new},
\ref{lem: verif ass Z}, 
\ref{lemma:euler:error:lp} and
\ref{lem: verif ass O} (see the Appendix below)
and Corollary \ref{cor: error exit times abstract_new}.

(b) To prove Theorem   \ref{thm:summary_for_sde T finite} we verify that condition 
\reff{eq: condition integrability for thm summary for sde} can be avoided when the time horizon is bounded. 
First we extend $\R^d$ to $\R^{d+1}$ equipped with the {Euclidean} metric and consider a function $\varrho\in C_b^2(\R)$ such that
\begin{enumerate}
\item $\varrho \le 0$ and $\varrho(0)=0$,
\item $\varrho$ is strictly increasing on $[-2L,0]$ and strictly decreasing on $[0,2L]$,
\item $\varrho \equiv -A$ on $[-2L,2L]^c$ for some $A>L$,
\item $D \varrho=1$ on $[-L-(r/2),-(r/2)]$ and $D \varrho=-1$ on $[(r/2),L+(r/2)]$.
\end{enumerate}
Note that our assumptions $0<r< 1/(4 L^3)$ and $L\ge 1$ guarantee the existence of such a $\varrho$.
Moreover, we can assume that $|\delta|_\infty\le L$ as $P$ and $\bar P$ take values in $[-L,L]$ only. 
We define the Lipschitz function $\delta^\#: \R^{d+1}\to \R$ by
\[ \delta^\# (x,y) := \delta(x) + \varrho(y) \]
and extend the open set $\Oc$ to an open set 
\[ \Oc^\# := \{ \delta^\# > 0 \} \subseteq \R^{d+1}. \]
By our construction we have that 
\begin{enumerate}
\item $\emptyset \not = \Oc^\# \subsetneq \overline  {\Oc^\#} \subsetneq \R^{d+1}$,
\item $\Oc^\# \subseteq \R^d \times [-2L,2L]$,
\item $\delta^\#$ is a distance function for $\Oc^\#$ in the sense of
      Assumption {\bf (P)}.
\end{enumerate}
Assume an auxiliary  one-dimensional Brownian motion $B=(B_t)_{t\ge 0}$ on a complete 
probability space $(\Omega',\Fc',\P')$ and define $\bar \Omega := \Omega\times \Omega'$ 
equipped with the completion $\bar \Fc$ of $\Fc \otimes \Fc'$ with respect to 
$\bar \P := \P \otimes \P'$. We extend the processes $W$, $B$, $X$ and $\bar X$ canonically 
to $\bar \Omega$ (where we keep the notation of the processes) and define the additional 
process $Y$ by $Y_t := B_{t\vee T}-B_T$. The right-continuous augmentation of the natural 
filtration of the $(d+1)$-dimensional Brownian motion
$(W,B)$ is denoted by $(\bar \Fc_t)_{t\ge 0}$. Therefore, letting  
$$
X^{\#}:=\left(\begin{array}{c}X\\Y\end{array}\right),           \hspace{.8em}
\bar X^{\#}:=\left(\begin{array}{c}\bar X\\Y\end{array}\right), \hspace{.8em}
                                                     \mbox{and} \hspace{.8em}
(P^{\#},\bar P^{\#}):=(\delta^{\#}(X^{\#}),\delta^{\#}(\bar X^{\#})),
$$ 
we obtain a setting that fulfills the assumptions of this paper. 
Now we check that 
$(X^{\#},\pi,P^{\#})$ and  $(\bar X^{\#},\bar \pi,\bar P^{\#})$ satisfy the conditions  \refAssZx,   \refAssOZx\  and \refassbornetaux\, 
with possibly modified parameters $(\kappa,L,r)$.
\smallskip

Assumption \refAssZx: For $(X^{\#},\pi,P^{\#})$  the assumption is trivial, the case 
$(\bar X^{\#},\bar \pi,\bar P^{\#})$ follows from the proof of Lemma \ref{lem: verif ass Z}.
\smallskip

Assumption \refAssOZx\ for $(X^{\#},\pi,P^{\#})$: The process $P^{\#}$ admits an It\^{o} decomposition 
\b*
dP^{\#}=b^{\#} dt + a^{\# \top} d\left(\begin{array}{c} W\\ B\end{array}\right), 
\e*
where $b^{\#}$ is uniformly bounded and 
$$
a^{\# \top}:=\left(\begin{array}{cc}D\delta(X) \sigma(X) & 0\\ (0,\ldots,0) & D \varrho(Y)\1_{[T,\infty)}\end{array}\right)
$$ is also bounded.
The condition  \refAssOZxii\ follows from the observation that $|\delta^{\#}|\le r/2$ implies either $|\delta|\le r$ or $r/2\le |\varrho|\le r/2+L$. 
\smallskip

Assumption \refAssOZx\ for $(\bar X^{\#},\bar \pi,\bar P^{\#})$: Similarly, using  \refAssOZxii\ for $r/2$, implies that 
$|\bar P_t| \le r$ and $|\bar X_t - \bar X_{\phi_t}| \le r/2$, or $r/2\le |\varrho|\le r/2+L$.
\smallskip

Assumption \refassbornetaux: It is sufficient to check the exit time of the process $Y$ from $[-2L,2L]$ computed on $\bar \pi$.
This follows by the arguments of the proof of Lemma \ref{lemma:both_conditions_L}.
\medskip

Finally we observe that $P_{t\wedge T} = P_{t\wedge T}^\#$, $\bar P_{t\wedge T} = \bar P_{t\wedge T}^\#$,
 $\theta_0 \wedge T = \theta_0^\# \wedge T$, and 
$\bar \theta_0^{\bar \pi} \wedge T = \bar \theta_0^{\bar \pi,\#} \wedge T$, where the quantities without $\#$ are taken with respect to
$(X,\bar X,\Oc)$ and the other ones for $(X^\#,\bar X^\#,\Oc^\#)$.
This implies that 
$\Esp{  |[\theta_{0}{  \wedge T}]-[\bar \theta^{\bar \pi}_{0}{  \wedge T}]|^p}\le  c_p \;   |\bar \pi|^{\frac12}$ and therefore 
$  \left(\Esp{  |X_{\theta_{0}{  \wedge T} }-\bar X_{ {\bar \theta^{\bar \pi}_{0}}{  \wedge T}}|^2}\right)^{\frac{ 1 }{2 }}
   \le c \;   |\bar \pi|^{\frac14}$. \qed


\appendix{}
\section{Appendix} 
 
\subsection{\bf Proof of Proposition \ref{proposition:reduction_bmo}.}
\underline{(L)$\Longleftrightarrow $(L'):}
 The condition \refassbornetaux\ obviously implies \refassbornetauxprime. Conversely, since 
$\Phi^{1,\pi}_{0}(\tau)\leq  \Econd{\tau}{\Phi^{1,\pi}_{0}(\phi^+_\tau)}+|\pi|$, where $|\pi|\le 1$,  and $\Phi^{1,\pi}_{0}(\phi^+_\tau)=0$ on $\{P_{\phi^+_\tau}\le 0\}$, the assumption 
\refassbornetauxprime\  implies that $\Phi^{1,\pi}_{0}(\tau)\leq L'+1$  for all $\tau \in \Tc$.
\vs2

\underline{(L)$\Longleftrightarrow $(L''):}
Indeed, \refassbornetaux\ implies \refassbornetauxsecond\ by Markov's inequality applied to the
level $c:=L/\alpha$ for a given $\alpha\in (0,1)$. Conversely, the fact  that $\theta^{ \pi }_{0}(\tau+ k c)=\theta^{ \pi }_{0}(\tau)$ on $\{\theta^{ \pi }_{0}(\tau)\geq \tau+ k c\}$ implies that 
$$
\P_\tau[\theta^{ \pi }_{0}(\tau)\geq \tau+ (k+1)c]=\E_\tau\big[\1_{\{\theta^{ \pi }_{0}(\tau)\geq \tau+ k c \}}\P_{\tau+ k c}[\theta^{ \pi }_{0}(\tau+ k c)\geq \tau+ (k+1)c]\big]. 
$$
Applying \refassbornetauxsecond\ inductively, allows us to conclude that the left-hand side above is controlled by $\alpha^{k+1}$. It follows that 
\begin{eqnarray*}
       \E_\tau[\theta^{ \pi }_{0}(\tau)-\tau]
&\leq& c +c \sum_{k\geq 0}\P_\tau[\theta^{ \pi }_{0}(\tau)\geq \tau+ (k+1)c]\\
&\leq& c +c \sum_{k\geq 0} \alpha ^{k+1}= c /(1-\alpha){ =: L}.
\end{eqnarray*} 
This proves that  \refassbornetauxsecond\ implies \refassbornetaux.  
\qed


\subsection{Proof of Lemma \ref{lem: error euler scheme up to stopping time with expo decay_new_new}}
(a) For $q>1$ we simply observe that 
\b*
      \left \| \Dist{X_{\vartheta},\bar X_{\vartheta}} \right \|_1
&\le& \sum_{k=0}^\infty\Esp{\sup_{t\in [k,k+1)} \Dist{X_{t},\bar X_{t}}
      \1_{\vartheta \in [k,k+1) }} \\
&\le&   \sum_{k=0}^\infty \Pro{\vartheta \in [k,k+1)}^\frac{q-1}{q}
                                \left \| \sup_{t\in [k,k+1)} \Dist{X_t,\bar X_t} \right \|_q \\
&\le&  \sum_{k=0}^\infty \Pro{\tau \ge k}^\frac{q-1}{q}
                                \left \| \sup_{t\in [k,k+1)} \Dist{X_t,\bar X_t} \right \|_q .
\e*
(b) follows immediately.
\ep


\subsection{\hspace*{-1em} Verification of the assumptions for the Euler scheme approximation}
\label{sec: proof thm euler}

All over this section, we work under the framework of Section \ref{sec: application to SDE}. We start with the condition \refAssZx\, 
that is - in a sense - independent from the set $\Oc$.
\vs2

\begin{lemma}\label{lem: verif ass Z} 
Under the Assumption \ref{ass: coeff mu sigma} the processes $X$ and $\bar X$ satisfy condition 
\refAssZx, where the function $\kappa:\R_+\times (0,\infty)\to \R_+$ depends at most on
$(L,d)$.
\end{lemma}
\proof For the process $X$ with the time-net $\pi=\R_+$, it  is trivially satisfied.
Let us fix $\tau \in \Tc$, $A\in \Fc_{\tau}$ of positive measure, $T>0$ and set 
$$ Y_{t}:=\bar X_{\tau + t T}\1_{A}
  \;\; \mbox{ for } t\in [0,1]. $$
For $2<p<\infty$ Assumption \ref{ass: coeff mu sigma} implies that 
$$
\E\left[|Y_{t}-Y_{s}|^{p}\right]\le c [T^p+T^\frac{p}{2}] |t-s|^{\frac{p}{2}} \P[A]
$$
for some $c=c(L,d,p)>0$ independent from the choice of $A\in \Fc_{\tau}$.
Fix $\alpha\in (0,\frac{1}{2}-\frac{1}{p})$. Then  it follows from the continuity of $Y$
and (the proof of) Kolmogorov's theorem in \cite[Theorem 2.1, p.26]{rev:yor:05} that
\b* 
      \E \left[\1_{A}\sup_{\tau\le t\le \tau+T} |\bar X_{t}-\bar X_{\bar \phi_{t}\vee \tau}|^{p}
      \right]
&\le& \E  \left[\1_{A}\sup_{|t-s|\le |\bar \pi|/T} | Y_{t}- Y_{s} |^{p}\right] \\
&\le& \frac{|\bar \pi |^{p\alpha}}{T^{p\alpha}} 
      \E  \left[\1_{A}\sup_{|t-s|\le |\bar \pi|/T; s,t\in D} 
      \frac{| Y_{t}- Y_{s} |^{p}}{(|\bar \pi |/T)^{p\alpha}}\right] \\
&\le& \frac{|\bar \pi |^{p\alpha}}{T^{p\alpha}} 
      \E  \left[\1_{A}\sup_{0\le s < t \le 1; s,t\in D} 
      \frac{| Y_{t}- Y_{s} |^{p}}{|t-s|^{p\alpha}}\right] \\
&\le& \frac{|\bar \pi |^{p\alpha}}{T^{p\alpha}} c' [T^p+T^\frac{p}{2}] \P[A],
\e*
where $c'=c'(c,p,\alpha)>0$ and $D\subseteq [0,1]$ are the dyadic points. Choosing  $p=6$ and $\alpha=1/6\in (0,\frac{1}{2}-\frac{1}{p})=(0,\frac{2}{6})$ gives 
$$   \Econd{\tau}{ \sup_{\tau\le t\le \tau+T} |\bar X_{t}-\bar X_{\bar \phi_{t}\vee \tau}|^{p}}
 \le c' [T^5+T^2] |\bar \pi | $$
and
$$
\P_{\tau}\left[\sup_{\tau\le t\le \tau+T} |\bar X_{t}-\bar X_{\bar \phi_{t}\vee \tau}|\ge \rho\right]\le \frac{c' [T^5+T^2]}{\rho^{p}} |\bar \pi|\;\;\mbox{  for $\rho>0$.}
$$
\qed
 
The next lemma is similar to \cite[Theorem A.1]{avikainen2007convergence}, which 
however involves a  $T^{2}$ term in the exponent, while we need a linear  term. 
It corresponds to the condition (b-ii) of Lemma \ref{lem: error euler scheme up to stopping time with expo decay_new_new}.
 
\begin{lemma}\label{lemma:euler:error:lp} 
If Assumption \ref{ass: coeff mu sigma} holds, then one has for all $4\le q<\infty$ {that} 
$$     \left \| \sup_{t\in [0,T]} |X_{t}-\bar X_{{t}}| \right \|_q
   \le Q(T,q) |\bar \pi |^\frac{1}{2},$$
where $Q(T,q) := c Q_q(T) e^{\alpha T}$ with $c>0$ depending at most on $(q,L,L_\mu,L_\sigma,d)$, a non-negative polynomial $Q_q$, and
$\alpha := d(6 L_\mu+3q L_\sigma^2)$.
\end{lemma}
\proof 
1. Let $2\le v < \infty$ and set $\Delta:=X-\bar X$. It follows from  It\^{o}'s Lemma that 
\b*
      |\Delta_{s}|^{2v}
& = & \int_{0}^{s} 2v |\Delta_{u}|^{2v-2} \Delta_{u}^{\top}d\Delta_{u} 
      +\sum_{i=1}^{d} \int_{0}^{s}  v |\Delta_{u}|^{2v-2} d\langle\Delta^{i}\rangle_{u} \\
&   & + \sum_{i,j=1}^{d} 2v(v-1) \int_{0}^{s} \Delta_{u}^{i}\Delta_{u}^{j} |\Delta_{u}|^{2v-4} 
      d\langle\Delta^{i},\Delta^{j}\rangle_{u}.
\e*
Under Assumption \ref{ass: coeff mu sigma} we obtain 
\b*
&   & \Esp{|\Delta_{s}|^{2v}}\\
&\le& \int_{0}^{s} 2v\Esp{ {|\Delta_{u}|^{2v-1}}|\mu(X_{u})
      -\mu(\bar X_{{\bar \phi}_{u}})|}du\\
&   & + \int_{0}^{s} v(1+2{d}(v-1))\Esp{ |\Delta_{u}|^{2v-2}
     |\sigma(X_{u})-\sigma(\bar X_{{\bar \phi}_{u}})|^{2}} du \\
&\le& A \int_{0}^{s} \Esp{ {|\Delta_{u}|^{2v-1}}
      \left( |\Delta_{u}| + |\bar X_{u}-\bar X_{{\bar \phi}_{u}}|\right) }du\\
&   & + B \int_{0}^{s} \Esp{ |\Delta_{u}|^{2v-2}
     | |\Delta_{u}| ^{2}+ |\bar X_{u}-\bar X_{{\bar \phi}_{u}}|^{2}  |} du \\
& = & [A+B]  \int_{0}^{s} \Esp{ |\Delta_{u}|^{2v}} du + 
      A \int_{0}^{s} \Esp{ {|\Delta_{u}|^{2v-1}}
      |\bar X_{u}-\bar X_{{\bar \phi}_{u}}| }du\\
&   & +  B \int_{0}^{s} \Esp{ |\Delta_{u}|^{2v-2}
       |\bar X_{u}-\bar X_{{\bar \phi}_{u}}|^{2}} du \\
\e*
for $A:= 2v L_\mu$ and $B:=2v(1+2{d}(v-1))L_\sigma^2 \le 6d v^2 L_\sigma^2$. Exploiting
$$ |\Delta_{u}|^{2v-1}
      |\bar X_{u}-\bar X_{{\bar \phi}_{u}}|
   \le \frac{2v-1}{2v} |\Delta_{u}|^{2v} + \frac{1}{2v} |\bar X_{u}-\bar X_{{\bar \phi}_{u}}|^{2v} $$
and 
$$ |\Delta_{u}|^{2v-2}
      |\bar X_{u}-\bar X_{{\bar \phi}_{u}}|^2
   \le \frac{v-1}{v} |\Delta_{u}|^{2v} + \frac{1}{v} |\bar X_{u}-\bar X_{{\bar \phi}_{u}}|^{2v} $$
we arrive at
\b*
&   & \Esp{|\Delta_{s}|^{2v}}\\
&\le& \left [ A+ B  + A \frac{2v-1}{2v} + B \frac{v-1}{v} \right ]
\int_{0}^{s} \Esp{ |\Delta_{u}|^{2v}} du \\
&   & + \left [ \frac{A}{2v} + \frac{B}{v} \right ]
          \int_{0}^{s} \Esp{|\bar X_{u}-\bar X_{{\bar \phi}_{u}}|^{2v} }du\\
&\le& 2 \left [ A+ B \right ]
\int_{0}^{s} \Esp{ |\Delta_{u}|^{2v}} du 
    + \left [ \frac{A}{2v} + \frac{B}{v} \right ]
          \int_{0}^{s} \Esp{|\bar X_{u}-\bar X_{{\bar \phi}_{u}}|^{2v} }du\\
&\le& 12 d [vL_\mu + v^2 L_\sigma^2]
       \int_{0}^{s} \Esp{ |\Delta_{u}|^{2v}} du \\
&   & + 6d[L_\mu+vL_\sigma^2]
      \int_{0}^{s} \Esp{|\bar X_{u}-\bar X_{{\bar \phi}_{u}}|^{2v} }du.
\e*
Exploiting
\b*
{\sup_{u\geq 0} }\; \Esp{ |\bar X_{u}-\bar X_{{\bar \phi}_{u}}|^{2v}}\le  c_v\;|\bar \pi|^{v}
\e*
for some constant $c_v=c(v,L)>0$, where we use the boundedness part of 
Assumption \ref{ass: coeff mu sigma}, we derive
$$  \Esp{|\Delta_{s}|^{2v}}
  \le 12 d [vL_\mu + v^2 L_\sigma^2]
       \int_{0}^{s} \Esp{ |\Delta_{u}|^{2v}} du +
      6d[L_\mu+vL_\sigma^2] s c_v  |\bar \pi|^v $$
and,  by Gronwall's Lemma,

$$  \Esp{|\Delta_{s}|^{2v}}
  \le  6d[L_\mu+vL_\sigma^2] c_v s e^{s 12 d[vL_\mu + v^2 L_\sigma^2]}
      |\bar \pi|^v. $$
2. Using the  It\^{o} decomposition of $|\Delta|^{2}$ and the Burkholder-Davis-Gundy and H\"older inequalities, we obtain (for another constant 
$c'_v=c'(v,L_\mu,L_\sigma)>0$) that
\b*
      \Esp{\sup_{0\le s\leq T}|\Delta_{s}|^{2v}}
&\le& c'_v \;[T^{v-1}+T^{v/2-1}] 
      \int_{0}^{T} \Esp{|\Delta_{u}|^{2v}+|\bar X_{u} - \bar X_{\bar \phi_{u}}|^{2v}
      }du \\
&\le& c'_v \;[T^v+T^{v/2}] |\bar \pi |^v 
      \big [ c_v + 6d[L_\mu+vL_\sigma^2] c_v e^{12 d T [vL_\mu + v^2 L_\sigma^2]}
            \big ] \\
&\le& c'_v c_v  \;[T^{v}+T^{v/2}]  
      \big [ 1 + 6d[L_\mu+vL_\sigma^2] \big ] e^{12 d T [vL_\mu + v^2 L_\sigma^2]} |\bar \pi |^v.
\e*
Consequently, for $q\ge 4$,
$$ \left \| \sup_{0\le s\leq T}|\Delta_{s}| \right \|_q
   \le C(q,L,L_\mu,L_\sigma,d) Q_q(T) e^{T d [6 L_\mu + 3 q  L_\sigma^2]} |\bar \pi |^\frac{1}{2}. 
   $$
\qed
\\
 
Now, we verify assumption \refAssOZx:

\begin{lemma}\label{lem: verif ass O} 
Let the Assumptions \ref{ass: coeff mu sigma}  and \ref{ass: non characteristic boundary condition} hold. Then $P$ and $\bar P$ satisfy the condition \refAssOZx\ for  $r>0$ small enough 
and $L\ge 1$ large enough, independently of $|\bar \pi|$.
\end{lemma} 
\proof 
First we apply It\^{o}'s Lemma to obtain that
$
d\bar P_{t}= \bar b_{t }  dt +  \bar a_{t}^{\top}  dW_{t}  
$
with 
$$
\bar b_{t }:=D\delta(\bar X_{t})\mu(\bar X_{{\bar \phi}_{t}})+\frac12{\rm Tr}[(\sigma\sigma^{\top})(\bar X_{{\bar \phi}_{t}}) D^{2}\delta(\bar X_{t})] 
\mbox{ and } \bar a_{t}^{\top}:= D\delta(\bar X_{t})\sigma(\bar X_{{\bar \phi}_{t}}).
$$
Up to an increase of $L$ in Assumption \refAssOZx\ (which potentially leads to a decrease
of $r$ to satisfy $0<r<1/(4L^3)$), condition  (i) is satisfied because 
$\delta,\mu,\sigma, D\delta, D^2 \delta$ are bounded.
 Since ${D\delta}$ is bounded by $L$ and $\sigma$ is $L$-Lipschitz, 
$$
    \left|D\delta(\bar X_{t})
    \sigma(\bar X_{{\bar \phi}_{t}})-D\delta(\bar X_{t})\sigma(\bar X_{t})\right|
\le L^{2}|\bar X_{{\bar \phi_t}}-\bar X_{t}|. 
$$
Consequently,
\b*
      |\bar a_t| 
&\ge& |D\delta(\bar X_{t})\sigma(\bar X_{t})| - 
       \left|D\delta(\bar X_{t})\sigma(\bar X_{{\bar \phi}_{t}})-
       D\delta(\bar X_{t})\sigma(\bar X_{t})\right| \\
&\ge& |D\delta(\bar X_{t})\sigma(\bar X_{t})| - L^2 |\bar X_{t}-\bar X_{{\bar \phi_t}}|.
\e*
For $|\bar P_t|\le r$, $|\bar X_{t}-\bar X_{{\bar \phi_t}}|\le r$ and
$0<r<1/(4L^3)$ this finally gives 
$|\bar a_t| \ge 1/L$ so that $\bar P$  satisfies \refAssOZxii.
The argument for $P$ is {analogous}.
\qed\\

We finally consider consider the Assumption {\bf (L)}.  
Conditions of type \eqref {eqn:exponential_tail_hitting_time} below can be found in 
\cite[Chapter 3]{frei:85}.

\begin{lemma}\label{lemma:both_conditions_L}
Let Assumption \ref{ass: coeff mu sigma} be satisfied and assume an $R>0$ and
  a non-increasing function $\varphi:[0,\infty)\to (0,\infty)$ with
$\lim_{T\to\infty} \varphi(T)=0$ such that
\begin{equation}\label{eqn:exponential_tail_hitting_time}
     \sup_{x\in \Oc}     \Pro{\theta_0^x(R) \ge T} \le  \varphi(T)
     \;\;\mbox{ for all } \;\; T\ge 0,
\end{equation}
where the open set $\Oc(R) := \Oc+ B_R$ ($B_R$ is the open ball centered at zero with radius $R>0$) 
satifies $\Oc(R)\subsetneq \overline {\Oc(R)} \subsetneq \R^d$
and $\theta_0^x(R):= \inf \{ t\ge 0 : X_t^x \not \in \Oc(R)\}$ for $x\in \Oc$ with  
$(X_t^x)_{t\ge 0}$ being the diffusion started in $x\in\R^d$.
Then there exist $\bar \eps\in (0,1]$ and a constant $K>0$ such that $|\bar \pi|\le \bar \eps$ implies
$$\Econd{\tau}{\bar \theta^{\bar\pi}_{0}(\tau)-\tau}+\Econd{\tau}{\theta_{0}(\tau)-\tau} \le K
         \;\;\mbox{ for all } \;\; \tau\in \Tc. $$
\end{lemma}
\proof 
      We only consider the estimate which involves $\bar \theta^{\bar\pi}_{0}(\tau)$
       (the other one follows directly from Proposition \ref{proposition:reduction_bmo}).
      By Proposition \ref{proposition:reduction_bmo} the case $\bar \theta^{\bar\pi}_{0}(\tau)$ can be reduced  
      to find $\alpha\in (0,1)$ and $c>0$ with
      \begin{equation}\label{eqn:Ltwo_is_sufficient}
            \Protau{\bar \theta^{\bar \pi}_{0}(\tau)-\tau \ge c} 
        \le \alpha \;\;  \mbox{for all}  \;\; \tau\in \Tc.
      \end{equation}
      Because for $c>1$ one has
      $$    \Protau{\bar \theta^{\bar \pi}_{0}(\tau)-\tau\ge c} 
        \le \Econd{\tau}{\Pcond{\bar\phi_\tau^+}{\bar \theta^{\bar \pi}_{0}(\bar\phi_\tau^+)-\bar\phi_\tau^+\ge c -1}}, $$
      it is sufficient to check \reff{eqn:Ltwo_is_sufficient} for $\tau\in \Tc^{\bar\pi}$.
      Given $\tau\in\Tc^{\bar\pi}$, we let
      $$\check X_t = x_0 + \int_0^t \check \mu_s ds + \int_0^t \check \sigma_s dW_s $$
      with $\check \mu_t := \1_{(0,\tau]}(t) \mu (\check X_{\bar \phi_t}) + 
                            \1_{(\tau,\infty)}(t) \mu (\check X_t)$ and with the corresponding definition for 
    $\check \sigma$. Let
      $\check \theta_0(\tau,R) := \inf \{ t\ge \tau : \check X_t \not \in\Oc(R) \}$.
      For $c\ge 2$ and $\tau\in \Tc^{\bar \pi}$ with $|\bar \pi|\le \bar \eps$, where $\bar \eps \in (0,1]$ is chosen at the
      end of the {proof},  we get from 
      {Lemma \ref{lemma:euler:error:lp}, applied to $q=4$ and some $T_0>0$,} for a set $A\in\Fc_\tau$ of positive 
      measure with $A\subseteq \{ \tau = t \} \cap \{ \bar X^{\bar \pi}_t \in \Oc \}$ (note that $\tau$ takes only 
      countable many values) that
      \b*
      &   & \Pro{A\cap \{ \bar \theta^{\bar \pi}_{0}(\tau)-\tau\ge c\} } \\
      & = & \Pro{A\cap \{ \bar \theta^{\bar \pi}_{0}(t)-t\ge c\} } \\
      &\le& \Pro{A\cap \{ \check \theta_0(t,R)-t \ge  c/2\} } \\
      &   & \hspace*{4em} + \Pro{A\cap\{ |\check X_{\check \theta_0(t,R)}-\bar X_{\check \theta_0(t,R)}|\ge R/2\}}\\
      &   & \hspace*{4em} + \Pro{A\cap \{ |\bar X_{\bar \phi^{+}_{\check \theta_0(t,R)}}-\bar X_{\check \theta_0(t,R)}|\ge R/2\}}\\
      &\le&   \Pro{A} \bigg [ \sup_{x\in \Oc} \Pro{\theta_0^x(R) \ge c/2}  \\
      &   & \hspace*{4em} 
            + \sup_{x\in \Oc} \sup_{|\tilde \pi|\le \bar\eps}
              \Pro{|X^x_{\theta_0^x(R)}- \bar X_{\theta_0^x(R)}^{x,\tilde \pi}|\ge R/2}
            +  \frac{c'(L)}{R^2} |\bar\pi| \bigg ] \\
      &\le& \Pro{A} \bigg [ \sup_{x\in \Oc} \Pro{\theta_0^x(R) \ge c/2}  {
            + \sup_{x\in \Oc} 
              \Pro{\theta_0^x(R) \ge T_0}} \\
      &   & { 
            + \sup_{x\in \Oc} \sup_{|\tilde \pi|\le \bar\eps}
              \Pro{|X^x_{\theta_0^x(R)\wedge T_0}- \bar X_{\theta_0^x(R)\wedge T_0}^{x,\tilde \pi}|\ge R/2}} 
            +  \frac{c'(L)}{R^2} |\bar\pi| \bigg ]  \\
      &\le& \Pro{A} \bigg [   { \varphi\left (\frac{c}{2}\right ) + \varphi(T_0) 
                  + \left ( \frac{2}{R} Q_\eqref{lemma:euler:error:lp}(T_0,4) \right )^4 \bar \eps^2
                  +  \frac{c'(L)}{R^2} \bar\eps } \bigg ]  ,
      \e*  
      where $\bar X^{x,\tilde \pi}$ is the Euler scheme for $X^x$ based on the net $\tilde \pi$. 
      First we choose $c\ge 2$ and $T_0>0$ large enough, then $\bar\eps$ small enough in order
      to arrange \reff{eqn:Ltwo_is_sufficient} for all $\tau\in \Tc^{\bar\pi}$ and some $\alpha \in (0,1)$.
 \qed
        

\end{document}